\documentclass[oneside,english]{amsart}
\usepackage[T1]{fontenc}
\usepackage[latin9]{inputenc}
\usepackage{amsthm,amsmath,amssymb,color,graphicx,babel,url}
\usepackage{float}
\usepackage{tikz}
\usepackage{graphicx}
\usepackage{mathtools}
\usepackage{makecell}

\usetikzlibrary{calc}

\newtheorem{theorem}{Theorem}
\newtheorem{lemma}[theorem]{Lemma}
\newtheorem{corollary}[theorem]{Corollary}
\newtheorem{question}[theorem]{Question}
\newtheorem{proposition}[theorem]{Proposition}
\theoremstyle{definition}
\newtheorem{definition}[theorem]{Definition}
\newtheorem{remark}[theorem]{Remark}
\newtheorem{example}[theorem]{Example}

\newcommand{\mQ}{{\mathcal Q}}
\newcommand{\mG}{{\mathcal G}}

\newcommand{\primozcomment}[1] {\color{magenta} ***  #1 *** \color{black}}
\newcommand{\micaelcomment}[1] {\color{blue} ***  #1 *** \color{black}}

\newcommand{\Aut}{\hbox{\rm Aut}}

\newcommand{\inv}{^{-1}}
\renewcommand{\mod}{\hbox{mod}\, }

\newcommand{\NN}{\mathbb{N}}
\newcommand{\ZZ}{\mathbb{Z}}

\newcommand{\Cay}{\hbox{\rm Cay}}
\newcommand{\Cov}{\hbox{\rm Cov}}
\newcommand{\GP}{\hbox{\rm GP}}

\newcommand{\fib}{\hbox{fib}}
\newcommand{\lcm}{\hbox{lcm}}

\newcommand{\beg}{\mathop{{\rm beg}}}
\newcommand{\term}{\mathop{{\rm end}}}
\renewcommand{\inv}{\mathop{{\rm inv}}}

\newcommand{\meo}{{\rm meo}}
\newcommand{\relmeo}{{\rm \eta}}

\newcommand{\rH}{{\rm H}}
\newcommand{\X}{{\rm X}}
\newcommand{\Y}{{\rm Y}}
\newcommand{\D}{{\rm D}}
\newcommand{\V}{{\rm V}}

\renewcommand{\lcm}{{\rm lcm}}

\newcommand{\Lad}{{\rm Moeb}}

\newcommand{\Prism}{{\rm Prism}}

\newcommand{\SDW}{{\rm SDW}}

\renewcommand{\mod}{\hbox{\rm{mod }}}

\makeatletter
\numberwithin{equation}{section}
\numberwithin{figure}{section}
\numberwithin{theorem}{section}

\makeatother

\usepackage{babel}

\numberwithin{equation}{section}
\numberwithin{figure}{section}

\title[Cubic vertex-transitive graphs]{Cubic vertex-transitive graphs admitting automorphisms of large order}
\author{Primo\v{z} Poto\v{c}nik}
\author{Micael Toledo}

\address{Primo\v{z} Poto\v{c}nik, Faculty of Mathematics and Physics, University of Ljubljana, Jadranska 21, SI-1000 Ljubljana, Slovenia.\newline
\indent Also affiliated with: Institute of Mathematics, Physics and Mechanics, Jadranska 19, SI-1000 Ljubljana, Slovenia.
}
\email{primoz.potocnik@fmf.uni-lj.si}

\address{Micael Toledo, D\'{e}partement de Math\'{e}matique, Universit\'{e} Libre de Bruxelles, C.P. 216 - Al\`{e}bre et Combinatoire, Boulevard du Triomphe, 1050, Brussels, Belgium}
\email{micael.alexi.toledo.roy@ulb.be}

\thanks{The authors gratefully acknowledge support of the Slovenian Research Agency: Core Programme P1-0294, Research Project J1-1691 and the Young Researcher Scholarship programme.}

\begin{document}

\nocite{*}

\begin{abstract}
A connected graph of order $n$ admitting a semiregular automorphism of order $n/k$ is called a $k$-multicirculant. Highly symmetric multicirculants of small valency have been extensively studied, and several classification results exist for cubic vertex- and arc-transitive multicirculants. In this paper we study the broader class of cubic vertex-transitive graphs of order $n$ admitting an automorphism of order $n/3$ or larger that may not be semiregular. In particular, we show that any such graph is either a $k$-multicirculant for some $k \leq 3$, or it belongs to an infinite family of graphs of girth $6$.

\end{abstract}

\maketitle

MSC2020: 05E18

\section{Introduction}

Studying the structure of the automorphism groups of highly symmetrical graphs is
one of the classical topics in the area of algebraic graph theory and a very important part of  it
aims at
proving upper bounds on the order of the automorphism groups in terms of a conveniently tame function of the order of the graph (see, for example, the classical work of Tutte~\cite{tutte} on cubic symmetric graphs).  Existence of such bounds often allows strong group theoretical tools to be applied.

In some  applications (such as a construction of a complete list of all
 graphs of fixed valence, bounded order and given symmetry type; see for example \cite{ConderFosterCensus,census}),
  a bound on the order of the automorphism group can be substituted with a weaker result where the order of individual automorphisms are bounded rather than the order of the automorphism group itself. 
It is well known that the order of the automorphism group of a connected vertex-transitive graph of valence $3$ ({\em cubic vertex-transitive graph}, for short) cannot be bounded by any subexponential function of the order (see, for example, \cite{census}). However, a recent result \cite[Theorem 1.6]{regorbs} shows that the order $o(g)$ of an individual automorphism $g$ of a cubic vertex-transitive graph other than $K_{3,3}$ equals the length $\ell(g)$ of the
longest orbit of the cyclic group $\langle g \rangle$, implying that $o(g)$
cannot exceed the order of the graph. In other words, if we define
$$
\relmeo(\Gamma) := \frac{ |\V(\Gamma)|}{\max \{ o(g) : g \in \Aut(\Gamma)\} },
$$
then 
$\relmeo(\Gamma) \ge 1$
 holds for every cubic vertex-transitive graph.
The aim of this paper is to investigate how sharp this bound is and under what additional assumptions it can be impoved.
More precisely, we 
obtain a complete classification of cubic vertex-transitive graphs $\Gamma$ 
for which 
$1\le \relmeo(\Gamma) \le 3$ 
holds and thus show that $\relmeo(\Gamma) > 3$ for every cubic vertex-transitive graph not
belonging to a known list of exceptional families
(see Theorem~\ref{the:main}).
 
 Before stating Theorem~\ref{the:main}, let us first introduce a few notions appearing in its statement.
A non-trivial automorphism $g$ of a graph is called {\em semiregular} provided that
the length of {\em every} vertex-orbit of $\langle g \rangle$ equals $o(g)$, or equivalently, when all the orbits of $\langle g \rangle$ have equal size.
 A graph admitting a semiregular automorphism with $k$ vertex-orbits is called a {\em $k$-multicirculant}; in addition,
 every graph on $n$ vertices is an $n$-multicirculant. Let 
 $$
 \kappa(\Gamma) := \min\{ k\, :\, \Gamma  \hbox{ is a } k\hbox{-multicirciulant} \}
$$
 and observe that
 $\relmeo(\Gamma)  \le \kappa(\Gamma) \le |\V(\Gamma)|$
 holds for every graph $\Gamma$.
What is more, the well-known {\em polycirculant conjecture}  \cite{Dragan} 
(which is known to be true for graphs of valence $3$ \cite{DraganSca}, as well as many other classes of graphs; see \cite{SpigaPCSurvey})
states that $\kappa(\Gamma) < |\V(\Gamma)|$ for every vertex-transitive graph $\Gamma$. 
 There are numerous classification results proved about cubic vertex-transitive $k$-multicirculants of different symmetry types  (see for example \cite{multi2,sym45,multi6,multi3,symtric,multi4})
and in particular, all cubic vertex-transitive $1$-multicirculants (also called {\em circulants}), $2$-multicirculant (also called {\em bicirculants})
and $3$-multicirculants (also called {\em triciculants})  are known \cite{bic,tricirc}.  

An interesting interplay between the parameters $\eta(\Gamma)$ and $\kappa(\Gamma)$ is considered in Section~\ref{sec:problem}.
In particular, the question whether the parameter $\kappa(\Gamma)$ can be bounded above in terms of $\eta(\Gamma)$ is discusses there.

Let us now introduce some families of graphs appearing in Theorem~\ref{the:main}.
For a positive integer $n\ge 3$, let $\Prism(n)$ denote the {\em prism} on $2n$ vertices, which can also be viewed as the generalised Petersen graph $\GP(n,1)$, and let
$\Lad(n)$, $n\ge 4$ even, be the {\em M\"obius ladder} with vertex-set $\ZZ_m$ and
edges of the form $\{x,x+s\}$ for $x\in \ZZ_n$, $s\in \{\pm 1, n/2\}$. Note that
 $\kappa(\Lad(n)) = 1$ for every $n$, while $\kappa(\Prism(n))$ is $1$ or $2$, depending on whether $n$ is odd or even, respectively.

Further, for a positive integer $n$ and distinct elements $i,j \in \ZZ_n\setminus \{0\}$, let $\rH(n,i,j)$ be the graph with vertex-set $\ZZ_n\times \ZZ_2$
and edges of the form $\{ (x,0), (x+s,1)\}$ for $x\in \ZZ_n$ and $s\in \{0,i,j\}$. Note that
the graphs $\rH(n,i,j)$ are bipartite bicirculants and are also known as
{\em cyclic Haar graphs}  \cite{bic}. Clearly, $\kappa(\rH(n,i,j)) \le 2$ and the values of $n,i,j$ for
which $\kappa(\rH(n,i,j)) = 1$ are characterised in Lemma~\ref{lem:haar}.

Let us now define two families of tricirculants, first introduced in \cite[Definitions 4.1 and 5.1]{tricirc}.
For an odd integer $k$, $k > 1$, let $\X(k)$ be the graph with $6k$ vertices, labelled 
$u_i$, $v_i$ and $w_i$ for $i\in \ZZ_{2k}$, and the edge-set being the union $E_1\cup E_2 \cup E_3 \cup E_4 \cup E_5$ where:
$E_1 = \{ \{u_i , u_{i+k}\} : i \in \ZZ_{2k}\}$,
$E_2 = \{ \{u_i,v_i\} : i \in \ZZ_{2k}\}$,
$E_3 = \{ \{u_i,w_i\} : i \in \ZZ_{2k}\}$,
$E_4 = \{ \{v_i,w_{i+1}\} : i \in \ZZ_{2k}\}$, and
$E_5 = \{ \{v_i,w_{i+r}\} : i \in \ZZ_{2k}\}$, where  $r=(k+3)/2$ if $k\equiv 1\> (\mod 4)$ 
and  $r=(k+3)/2 + k$ if $k\equiv 3\> (\mod 4)$.

Further, for an odd integer $k$, $k > 1$, let $\Y(k)$ be the graph with $6k$ vertices, labelled 
$u_i$, $v_i$ and $w_i$ for $i\in \ZZ_{2k}$, and the edge-set being the union $E_1\cup E_2 \cup E_3 \cup E_4 \cup E_5$ where:
$E_1 = \{ \{u_i , u_{i+1}\} : i \in \ZZ_{2k}\}$,
$E_2 = \{ \{u_i,v_i\} : i \in \ZZ_{2k}\}$,
$E_3 = \{ \{v_i,w_i\} : i \in \ZZ_{2k}\}$,
$E_4 = \{ \{v_i,w_{i+2}\} : i \in \ZZ_{2k}\}$, and
$E_5 = \{ \{w_i,w_{i+k}\} : i \in \ZZ_{2k}\}$.

By \cite[Theorem 4.3 and Theorem 5.3]{tricirc}, 
 $\kappa(\X(k))=\kappa(\Y(k)) =3$ for $k \equiv 3 \> (\mod 6)$ while
  $\kappa(\X(k))=\kappa(\Y(k)) =2$ for $k \equiv 0 \> (\mod 6)$.
 More facts about graphs $\X(k)$ and $\Y(k)$ can be found in 
 \cite[Sections 4 and 5]{tricirc}.

 For positive integers $m, t \ge 3$, let $\SDW(m,t)$ be
the cubic graph with the vertex-set $\ZZ_m\times \ZZ_t \times \ZZ_2$ and edges
of the form $\{(x,i,0),(x,i\pm 1,1)\}$ and $\{(x,i,1),(x+1,i,0)\}$ for all $x\in \ZZ_m$ and $i\in \ZZ_t$.
(The graphs $\SDW(m,t)$ with $t=3$ are also known as {\em split depleted wreath graphs}.)
Note that a graph $\SDW(m,t)$ is a $2t$-multicirculant, as witnessed by the semiregular automorphism of order $m$, mapping a vertex $(x,i,j)$ to the vertex $(x+1,i,j)$ for every $x,i,j$.
It can be easily seen that $\kappa(\SDW(m,3)) \le 6$ and that
$\relmeo(\SDW(m,3)) \le 3$ unless $m$ is divisible by $6$ (in which case
$\relmeo(\SDW(m,3)) = 6$); see Lemma~\ref{lem:deltanomulti} for more details.

Finally, Tutte's $8$-cage is the unique cubic arc-transitive graph on $30$ vertices, appearing under the name {\tt CubicVTgraph}$(30,8)$ in \cite{census}, while the truncated tetrahedron is the graph on $12$ vertices obtained by the geometric truncation of the skeleton of the tetrahedron (it appears under the name {\tt CubicVTgraph}$(12,2)$ in \cite{census}).

We can now state the main result of this paper.

\begin{theorem}
\label{the:main}
Let $\Gamma$ be a finite simple connected vertex-transitive graph of valence $3$ and order $n$. Then $\Gamma$ admits an automorphism of order at least $\frac{n}{3}$
(or equivalently, $\relmeo(\Gamma) \le 3$) if and only if one of the following happens:
\smallskip
\begin{enumerate}
\item $\kappa(\Gamma) = 1$  and $\Gamma$ is isomorphic to
\begin{enumerate}
 \item the prism $\Prism(m)$ where $n=2m$ with $m\ge 3$, $m$ odd; or
 \item the M\"obius ladder $\Lad(n)$ with $n\ge 4$;
\end{enumerate}
\smallskip
 \item $\kappa(\Gamma) =2$ and $\Gamma$ is isomorphic to
\begin{enumerate}
 \item the prism $\Prism(m)$ where $n=2m$ with $m\ge 4$, $m$ even;
 \item a generalised Petersen graph $\GP(m,r)$ where $n=2m$, $m\ge 5$, $2\le r < m/2$, and
   \begin{itemize}
    \item $m \geq 3$, $r^2 \equiv \pm 1$ $(\mod m)$, or
    \item $m = 10$ and $r=2$;
   \end{itemize} 
 \item the cyclic Haar graphs $\rH(m;r,s)$ where $n=2m$ with 
  $m \geq 3$, $1 \le r, s \le m-1$, $r\not = s$, $r$ divides $m$, $\gcd(r,s)=1$,
  such that $m$ is even or $m$ is odd and $\{r,s\}$ is neither $\{1,m-1\}$ nor $\{1,2\}$.
\end{enumerate}
\smallskip
\item $\kappa(\Gamma) = 3$ and $\Gamma$ is isomorphic to one of the following graphs:
\begin{enumerate}
\item $\X(k)$ with $k\equiv 3\> (\mod 6)$;
\item $\Y(k)$ with $k\equiv 3\> (\mod 6)$;
\item Tutte's 8-cage where $n=30$;
\item the truncated tetrahedron where $n=12$.
\end{enumerate}
\smallskip
\item $\kappa(\Gamma) = 6$ and $\Gamma \cong \SDW(m,3)$ with $m\equiv 3$ $(\mod 6)$,  
$m\ge 9$,
 $n=6m$.
\end{enumerate}
\end{theorem}

\begin{remark}
If $\Gamma$ is a graph appearing in one of items (1)--(3), then $\kappa(\Gamma) =  \relmeo(\Gamma)$ except when $\Gamma$ is isomorphic to 
\begin{itemize}
\item the cube graph ${\rm Q}_3 \cong \Prism(4)$ where $\relmeo(\Gamma) = \frac{4}{3}$;
\item the Petersen graph $\GP(5,2)$ where $\relmeo(\Gamma) = \frac{5}{3}$;
\item the Heawood graph $\rH(7,1,3)$ where $\relmeo(\Gamma) = \frac{7}{4}$;
\item the M\"{o}ebius-Kantor graph $\GP(8,3)$ where $\relmeo(\Gamma) = \frac{4}{3}$;
\item the Pappus graph $\Y(3) \cong \SDW(3,3)$ where $\relmeo(\Gamma) = \frac{3}{2}$.

\end{itemize}
If $\Gamma$ is one of the graphs in item (4) then $\relmeo(\Gamma) = 3$. The relation between the functions $\relmeo$ and $\kappa$ is further discussed in Section \ref{sec:problem}.
\end{remark}

 The proof of the above theorem is inevitably rather technical since a certain amount of case-by-case analysis cannot be avoided. However, we have tried to use as many theoretical tools
as possible in order to shorten and organise the arguments into self contained parts.

 Our proof relies on two crucial ideas. The first idea is to classify  the cubic vertex-transitive graphs $\Gamma$ admitting an automorphism $g$ of order at least $|\V(\Gamma)|/3$ in terms of the quotient graphs $\Gamma/\langle g \rangle$. For this approach to be practical, we need to
 store enough information that will allow us to reconstruct the graphs $\Gamma$ from
 their quotients. Here, a recently developed theory of generalised cyclic covering projections
 \cite{cycliccovers}, briefly summarised in Section \ref{sec:covers}, provided the needed theoretical background.
The second crucial fact making this approach feasible follows from the results of \cite{regorbs},
from which one can deduce that a cyclic group of automorphisms
of a cubic vertex-transitive graph $\Gamma$ with order at least $\frac{|\V(\Gamma)|}{3}$ can 
have at most $5$ orbits on $\V(\Gamma)$ (see Lemma~\ref{lem:bound5}).
 In particular, there is only a finite number of possible quotients $\Gamma/\langle g \rangle$
with $o(g) \ge \frac{|\V(\Gamma)|}{3}$. The rest of the proof is then a careful
case-by-case analysis of these possible quotient graphs. The strategy of proving Theorem~\ref{the:main} is laid out in more details in Section~\ref{sec:strategy}.

\section{Overview and basic definitions}
\label{sec:2}

The following paragraphs, in which we give some basic formal definitions and we outline the proof of Theorem \ref{the:main}, serve as a more detailed summary of the contents of this paper.

\subsection{Graphs}
\label{ssec:graphs}
We would first like to stress that all the graphs in this paper are finite.
Even though we are primarily interested in simple graphs (that can be defined as a finite set
of vertices together with an irreflexive symmetric relation on it),
it will be very convenient for us to adopt a more general definition of a graph 
that has become standard when quotients and covers of
graphs are considered (see, for example, \cite{covers}).

For us, a {\em graph} is an ordered $4$-tuple $(V,D; \beg,\inv)$ where
$D$ and $V \neq \emptyset$ are disjoint finite sets of {\em darts}
and {\em vertices}, respectively, $\beg\colon D \to V$ is a mapping
which assigns to each dart $x$ its {\em initial vertex}
$\beg\,x$, and $\inv$ is an involutory permutation of $D$ which interchanges
every dart $x$ with its {\em inverse dart}, also denoted by $x^{-1}$.
The final vertex of a dart $x$ is $\beg x^{-1}$ and is denoted $\term x$.
The {\it neighbourhood} of a vertex $v$ is defined as the set
of darts that have $v$ for its initial vertex and the
{\it valence} of $v$ is the cardinality of the neighbourhood. 
If $\Gamma$ is a graph we write $\V(\Gamma)$ and $\D(\Gamma)$ to denote the vertex- and dart-set of $\Gamma$, respectively. Furthermore, we may write $\beg_{\Gamma}$ and $\inv_{\Gamma}$, with a subscript, to indicate the beginning and inverse functions of $\Gamma$, to avoid confusion when more than one graph is involved. We will generally omit the subscript if there is no possibility of ambiguity.

The orbits of $\inv$ are called {\em edges}.
The edge $\{x,x^{-1}\}$ containing a dart $x$ is called a {\em semiedge} if $x^{-1} = x$,
a {\em loop} if $x^{-1} \neq x$ while
$\beg\,(x^{-1}) = \beg\,x$,
and  is called  a {\em link} otherwise.
The {\em endvertices of an edge} are the initial vertices of the darts contained in the edge. If $\{u,v\}$ is the set of the endvertices of an edge, then we say that $u$ and $v$ are adjacent and write $u \sim v$.
Two darts $x$ and $y$ are {\em parallel} if $\beg x = \beg y$ and $\beg x^{-1} = \beg y^{-1}$. Two edges are parallel if they have the same endvertices. 
When we present a graph as a drawing, the links are drawn in the usual way as a line
between the points representing its endvertices, a loops is drawn as a closed curve at its unique endvertex and a semiedge is drawn as a segment attached to its unique endvertex. 

A graph without loops, semiedges and pairs of parallel edges is {\em simple}. Note that a simple graph is completely determined by its vertex-set and the adjacency relation, and conversely, given
a set $V$ and an irreflexive symmetric relation $\sim$ on $V$, we can define a graph 
$(V,D; \beg,\inv)$  by letting $D := \{ (u,v) \colon u,v \in \V(\Gamma), u\sim v\}$,
 $\beg(u,v) = u$ and $\inv(u,v) = (v,u)$. A dart in a simple graph is traditionally called an {\em arc},
 so we will use these two terms interchangeably. 
 In this paper, a {\em cubic graph} will always stand for a connected simple graph in which every vertex has valence $3$.
 
Notions such as morphism, isomorphism and automorphism of graphs are obvious
generalisations of those defined in the traditional setting and precise definitions can be found in \cite{covers,elabcovers}. In particular, an automorphism of a graph $(V,D; \beg,\inv)$
is a permutation $g$ of $V\cup D$ preserving each of $V$ and $D$ such that
$\inv(x^g) = \inv(x)^g$ and $\beg(x^g) = \beg(x)^g$ for every $x\in D$. Note that
when the graph is simple, an automorphism is uniquely defined by its adjacency preserving action on the vertex-set. We shall thus often consider view automorphism of simple graphs in the usual way, that is, as adjacency preserving permutations of the vertex-set.

\subsection{Labelled quotients of graphs}
Let $\Gamma:=(V,D; \beg,\inv)$ be a graph admitting a cyclic group of automorphisms $G \leq \Aut(\Gamma)$. For $v \in V$, let $v^G$ denote the $G$-orbit of $v$ and let $V/G$ be the set of all $G$-orbits of vertices of $\Gamma$. Similarly, let $x^G$ be the $G$-orbit of $x \in D$ and let $D/G$ be the set of all $G$-orbits on darts. We define the $G$-quotient of $\Gamma$ as the graph $\Gamma/G=(V/G,D/G,\beg',\inv')$ where $\beg' x^G = (\beg x)^G$ and $\inv' x^G = (\inv x)^G$ for all $x \in \D(\Gamma)$. 

\begin{figure}[h!]
\centering
\includegraphics[width=0.6\textwidth]{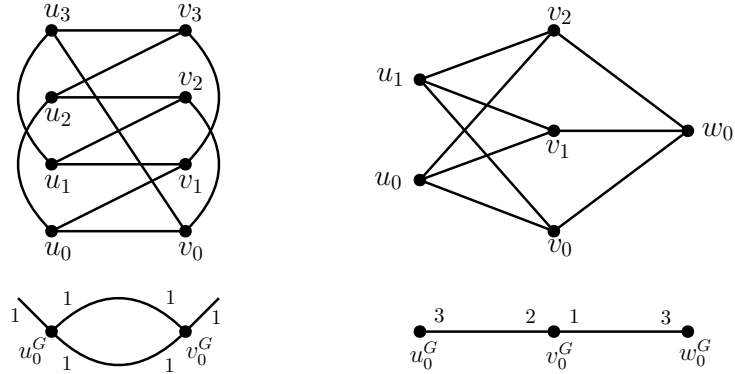}
\caption{The M\"{o}bius ladder on $8$ vertices (top left) and $K_{3,3}$ (top right) above their labelled quotients by a cyclic group of automorphisms $G = \langle g \rangle$. For the case on the left, $g$ is the automorphisms given by the permutation $(u_0u_1u_2u_3)(v_0v_1v_2v_3)$; for the graph on the right $g$ is given by $(u_0u_1)(v_0v_1v_2)$. In the graph on the lower right, the edge joining $u_0^G$ and $v_0^G$ has a label next to each of its endvertices: $3$ next to $u_0^G$, and $2$ next to $v_0^G$. These are the labels $\lambda_G(x)$ and $\lambda_G(x^{-1})$ where $x$ is the dart pointing from $u_0$ to $v_0$.
 This indicates that every vertex in the fibre of $u_0^G$ (in the graph above) has $3$ neighbours in the fibre of $v_0^G$ while every vertex in the fibre of $v_0^G$ has $2$ neighbours in the fibre of $u_0^G$.   }
\label{fig:cocientes}
\end{figure}

Let $x \in \D(\Gamma)$ be a dart and let $u = \beg x$. Then, $x^G$ is a dart of $\Gamma / G$ with initial vertex $u^G$. Let $\lambda_G(x^G)$ denote the number of darts of $\Gamma$ in the orbit $x^G$ that begin at any fixed vertex in $u^G$ (note that this is independent of which vertex of $u^G$ we choose and that $\lambda_G(x^G) = |x^{G_u}|$). Then $\lambda_G \colon \D(\Gamma/G) \to \NN$ is a well-defined function and the pair $(\Gamma/G , \lambda_G)$ is called a labelled quotient (see Figure \ref{fig:cocientes}). 

\subsection{Strategy}
\label{sec:strategy}
In what follows, we describe our plan to prove Theorem \ref{the:main}. 
We first consult the census of all cubic vertex-transitive graphs \cite{census} and check  that the theorem holds for all the graphs on at most $20$ vertices (this can be easily done with a help of computer and a computer algebra system such as Sage~\cite{sage}).

We may thus concentrate on the class
 $\mG$ of all cubic vertex-transitive graphs on at leats $n > 20$ admitting a cyclic group $G \leq \Aut(\Gamma)$ of order at least $\frac{n}{3}$. Let $\mQ$ be the set of labelled quotients $(\Gamma/G, \lambda_G)$ where $\Gamma \in \mG$ and $G \leq \Aut(\Gamma)$ is cyclic of order at least $\frac{|\V(\Gamma)|}{3}$ (by \cite[Theorem 4.7]{regorbs}, the order of $G$ equals the order of the largest $G$-orbit on vertices). 

Observe that if $(\Gamma/G, \lambda_G) \in \mQ$, then the graph $\Gamma/G$ is connected but may admit parallel edges, loops or semi-edges.
Since the valence of every vertex in $\Gamma$ is $3$, it follows that the vertices
in $\Gamma/G$ have valence at most $3$ and that
 $\lambda_G(x) \leq 3$ for every dart $x$ of $\Gamma/G$. 
Moreover, a labelled graph in $\mQ$ can have at most $5$ vertices (see Lemma \ref{lem:bound5}).
 This shows that the set $\mQ$ is a subset of the set $\mQ^\circ$ of all connected subcubic labelled graphs $(Q,\lambda)$ on at most $5$ vertices with $\lambda(x) \le 3$. The set $\mQ^\circ$ is clearly
finite, but still consists of an inconveniently large number of labelled graphs.
 
To determine which of the labelled graph in $\mQ^\circ$ 
indeed arise as quotients of the cubic graphs in $\mG$ by an appropriate cyclic group $G$ (that is, which of them belong to $\mQ$), we rely on the concept of a cyclic generalised voltage graph, first introduced in \cite{cycliccovers}, and the associated generalised covering graph construction. 
Loosely speaking, this construction takes a labelled graph $(\Delta,\lambda)$, together with an additional information, called {\em voltage assignment} $\zeta$, as an input and constructs a connected graph $\Gamma$ (called a {\em cover}, for short) having $(\Delta,\lambda)$ as a labelled quotient. It was proved in \cite{cycliccovers}, that as the voltage assignment $\zeta$ varies  this procedure yields all possible connected graphs having $(\Delta,\lambda)$ as a quotient by a cyclic group. 

In Section \ref{sec:VT}, we use the results proved in \cite{cycliccovers} (and summarised in Section~\ref{sec:covers}) to find a set of necessary conditions (see Theorem~\ref{the:diagram}) for a labelled graph in $\mQ^\circ$ to admit a connected vertex-transitive generalised cyclic cover belonging to $\mG$. In Section \ref{sec:obstructions} we determine, by means of forbidden labelled subgraphs, further necessary conditions for an element of $\mQ^\circ$ to admit a cyclic  generalised cover in $\mG$ (see Theorem~\ref{the:artefacts}). The set of conditions given in Theorems \ref{the:diagram} and \ref{the:artefacts} is restrictive enough to allow us to compute, via a brute-force algorithm, the subset $\mQ^* \subseteq \mQ^\circ$ of labelled graphs satisfying them. There are $20$ such graphs.

\begin{figure}[h!]
\centering
\includegraphics[width=0.75\textwidth]{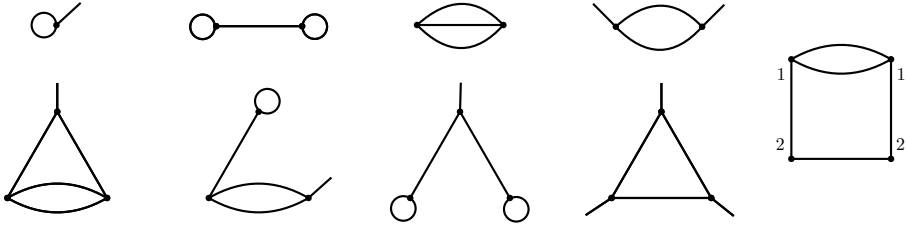}
\caption{The nine labelled graphs in $\mQ$, where darts with no label have $\lambda$-value $1$. With the exception of the right-most graph, these graphs correspond to the $8$ possible quotients of a $k$-multicirculant graph by a $k$-multicirculant automorphism, with $k \in \{1,2,3\}$. }
\label{fig:Q}
\end{figure}

In Section \ref{sec:deltas}, we analyse the $20$ labelled graphs of $\mQ^*$ in detail, and show that precisely nine of them admit vertex-transitive cubic cyclic generalised covers belonging to $\mG$ (see Figure~\ref{fig:Q}). These are the nine elements of the set $\mQ$. However, due to some overlap in the families of covering graphs, only seven elements of $\mQ$ are necessary to reconstruct $\mG$. 
Finally in Section \ref{sec:last}, we characterise the elements of $\mG$ and complete the proof of Theorem \ref{the:main}. 

\section{Covers}
\label{sec:covers}

We now formally introduce the concept of a {\em cyclic generalised voltage graph}, which generalises voltage graphs (in the sense of \cite{grosstuc}) for cyclic voltage groups, and is a special case of the wider class of generalised voltage graphs introduced in \cite{MPgc}. The definitions and results in this section are mostly taken from \cite{cycliccovers},
where cyclic generalised voltage graphs were first defined.
 Each cyclic generalised voltage graph gives rise to a unique {\em generalised covering graph} (called {\em covering graph} for simplicity). By the end of the section, we characterise those cyclic generalised voltage graphs whose covering graphs are cubic (that is,  connected, finite, simple $3$-valent graphs).

%
%
%
%
%
%
%
%
%
%
%

\begin{definition}
\label{def:cgvg}
Let $\Delta$ be a finite connected graph and let $\lambda\colon \D(\Delta) \to \NN$, $\iota\colon \V(\Delta) \to \NN$ and $\zeta\colon \D(\Delta) \to \ZZ$ be functions such that
\begin{eqnarray}
\label{eq:ratio}\lambda(x)\iota(\beg x) &=& \lambda(x^{-1})\iota(\beg x^{-1}),\\
\label{eq:invvolt}\zeta(x^{-1}) &\equiv& -\zeta(x) \quad (\mod \lambda(x)\iota(\beg x))
\end{eqnarray}
for every dart $x\in \D(\Delta)$. Then we say that
the quadruple $(\Delta,\lambda,\iota,\zeta)$ is a {\em cyclic generalised voltage graph}, and we call the functions $\lambda$, $\iota$ and $\zeta$ a {\em labelling}, an {\em index function} and a {\em voltage assignment}, respectively.
\end{definition}

\begin{definition}
Let $(\Delta,\lambda,\iota,\zeta)$ be a cyclic generalised voltage graph. The cover of $(\Delta,\lambda,\iota,\zeta)$, denoted $\Cov(\Delta,\lambda,\iota,\zeta)$ is the graph $\Gamma := (V',D',\beg',\inv')$ where:

\begin{enumerate}
\item $V' = \{(v,i) \colon v \in \V(\Delta), i \in \{0,\ldots,\iota(v)-1\}\}$;
\item $D' = \{(x,i) \colon x \in \D(\Delta), i \in \{0,\ldots,\lambda(x)\iota(\beg x)-1\}\}$;
\item $\beg' (x,i) = (\beg x,i)$;
\item $(x,i)^{-1} = (x^{-1},i + \zeta(x))$.
\end{enumerate}
\end{definition}

\begin{remark}
Since the second coordinate of a vertex $(v,i)$ or a dart $(x,i)$ in the definition above is an element of $\{0,\ldots,\iota(v)-1\}$ or $\{0,\ldots,\lambda(x)\iota(\beg x)-1\}$, respectively, any operation on the second coordinate is to be computed modulo $\iota(v)$ or $\lambda(x)\iota(\beg x)$ accordingly. In particular, on the right-hand side of equality (4), the sum $i + \zeta(x)$ is to be computed modulo $\lambda(x)\iota(\beg x)$.
\end{remark}

Let $(\Delta,\lambda,\iota,\zeta)$ be a cyclic generalised voltage graph and set $\Gamma = \Cov(\Delta,\lambda,\iota,\zeta)$. For the sake of simplicity, we will write $v_i$ instead of $(v,i)$ for a vertex of $\Gamma$, and $x_i$ instead of $(x,i)$ for a dart of $\Gamma$. The natural projection $\pi: \Gamma \to \Delta$ that maps every $v_i \in \V(\Gamma)$ to $v$ and every $x_i \in \D(\Gamma)$ to $x$ is a graph epimorphism. For each vertex $v \in \V(\Delta)$ we call the set $\pi^{-1}(v)=\{v_i \colon i \in \{1,\ldots,\iota(v)\}\}$ the {\em fibre} of $v$ and we denote it $\fib(v)$. Similarly, the {\em fibre} of a dart $x \in \D(\Delta)$ is $\fib(x)=\pi^{-1}(x)=\{x_i \colon i \in \{0,\ldots,\lambda(x)\iota(\beg x)-1\}\}$. 

Let $n= \lcm\{\lambda(x)\iota(\beg x) \colon x \in \D(\Delta)\}$ and observe that the group $\ZZ_n$ acts on $\D(\Gamma) \cup \V(\Gamma)$ by the rule $(z_i)^g=z_{i+g}$  for all $z_i \in \D(\Gamma) \cup \V(\Gamma)$ and $g \in \ZZ_n$ (recall that the index $i+g$ is computed modulo $\iota(z)$, if $z \in \V(\Gamma)$, or modulo $\lambda(z)\iota(\beg z)$ if $z \in \D(\Gamma)$). The permutation induced by each $g \in \ZZ_n$ is an automorphism of $\Gamma$ and the orbit of a dart (or a vertex) $z_i$ under this action is precisely $\fib(z)$. Furthermore, the action of $\ZZ_n$ is faithful and hence there is an embedding $\phi:\ZZ_n \to \Aut(\Gamma)$ (see \cite[Lemma 6.3]{cycliccovers}).
The image $\phi(\ZZ_n)$ is a cyclic group of order $n$ generated by the automorphism of $\Gamma$ that maps every vertex $v_i \in \V(\Gamma)$ to $v_{i+1}$ and every dart $x_i \in \D(\Gamma)$ to $x_{i+1}$. We call this automorphism the {\em canonical covering transformation} of $\Cov(\Delta,\lambda,\iota,\zeta)$ and we denote it by $\rho$.

 In short, the covering graph $\Gamma$ admits a cyclic group of automorphisms of order $n$ whose orbits on vertices and darts are precisely the fibres of vertices and darts of $(\Delta,\lambda,\iota,\zeta)$. 

Conversely, in view of \cite[Theorem 5.3]{cycliccovers} and \cite[Theorem 6.2]{cycliccovers}
every graph admitting a cyclic subgroup of automorphism of order $n$ is the cover of some cyclic generalised voltage graph $(\Delta,\lambda,\iota,\zeta)$,
and by \cite[Lemma 6.1]{cycliccovers}, it follows that $\lcm\{\lambda(x)\iota(\beg x) \colon x \in \D(\Delta)\} = n$. Moreover, by \cite[Theorem 6.6]{cycliccovers}  we can always assume that $\zeta(x) = 0$ for all darts $x$ lying on a prescribed spanning tree $\mathcal{T}$ of $\Delta$. A voltage assignment satisfying this condition is said to be {\em $\mathcal{T}$-normalised}.
Let us summarise these observations in the following theorem:

\begin{theorem}
\label{theo:cover}
A graph $\Gamma$ admits a cyclic subgroup of automorphisms $G \leq \Aut(\Gamma)$ of order $n$ if and only if $\Gamma \cong \Cov(\Gamma/G, \lambda,\iota,\zeta)$ for some functions $\lambda$, $\iota$ and $\zeta$ where $\zeta$ is $T$-normalised for a spanning tree $T$ of $\Gamma /G$ and $n= \lcm\{\lambda(x)\iota(\beg x) \colon x \in \D(\Delta)\}$.
\end{theorem}

\begin{figure}[h!]
\centering
\includegraphics[width=0.4\textwidth]{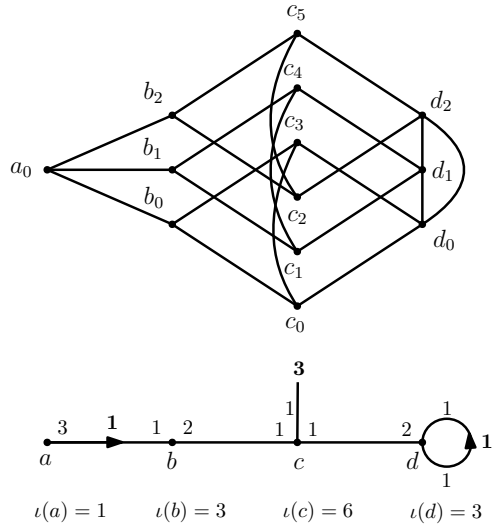}
\caption{A cyclic generalised voltage graph (bottom) along with its cyclic generalised cover (top).}
\label{fig:firstexample}
\end{figure}

\begin{example}
Consider the cyclic generalised voltage graph depicted in the bottom of Figure \ref{fig:firstexample}. The function $\lambda$ is given in the figure as follows: each link has two numbers next to each of its endvertices, corresponding to the $\lambda$-values of each dart underlying this link; for instance, if we let $x$ be the dart beginning at $a$ and ending at $b$, then $\lambda(x) = 3$ and $\lambda(x^{-1})=1$. The semi-edge at $c$ has a single label $1$ and the loop at $d$ has two labels, both equal to $1$, corresponding to the $\lambda$-values of the darts underlying it. The voltage assignment is also given in the figure. The link joining $a$ with $b$ has an arrowhead directed from $a$ to $b$ with a $1$ in boldface written above it. This indicates that the dart beginning at $a$ and ending at $b$ has voltage $1$, while its inverse $x^{-1}$ has voltage $-1$. The semi-edge at $c$ has voltage $3$, as indicated by the boldface number above it. One dart underlying the loop at $d$ has voltage $1$ while its inverse has voltage $-1$. The values of $\iota$ are written below the graph. The graph at the top of Figure \ref{fig:firstexample} is then the cover of the cyclic generalised voltage graph.

\end{example}

\begin{remark}
Let $(\Delta,\lambda,\iota,\zeta)$ be a cyclic generalised voltage graph and let $\Gamma = \Cov(\Delta,\lambda,\iota,\zeta)$. Let $x \in \D(\Delta)$ and $u = \beg x$. Then the following is straightforward from the definition of a cyclic generalised voltage graph:
\begin{enumerate}
\item $|\fib(u)| = \iota(u)$;
\item $|\fib(x)| = \lambda(x)\iota(u)$;
\item For every $u_i \in \fib(u)$, there are exactly $\lambda(x)$ darts in $\fib(x)$ that begin at $u_i$;
\item For every $u_i \in \fib(u)$, $\deg(u_i) =\deg_{\lambda}(v)$, where
\end{enumerate}
\end{remark}
\begin{align}
\deg_{\lambda}(v) = \sum\limits_{x \in \Delta(v)} \lambda(x).
\end{align} 

Since cubic graphs are the main object of study of this paper, we would like to focus precisely on those cyclic generalised voltage graphs whose covers are cubic graphs. 
Hence the following definition.

\begin{definition}
\label{def:ccv}
A cyclic generalised voltage graph $(\Delta,\lambda,\iota,\zeta)$ is called a $ccv$-graph whenever $\Cov(\Delta,\lambda,\iota,\zeta)$ is a cubic graph. 
\end{definition}

The following characterization of $ccv$-graphs is a consequence of  \cite[Theorems 6.8 and 6.9]{cycliccovers}.

\begin{lemma}
\label{lem:ccv}
Let $(\Delta,\lambda,\iota,\zeta)$ be a cyclic generalised voltage graph where $\zeta$ is $\mathcal{T}$-normalised for some spanning tree $\mathcal{T}$ of $\Delta$. Let $A = \{\zeta(x) \colon x \in \D(\Delta)\}$ and $B = \{\iota(v) \colon x \in \V(\Delta)\}$. Then $(\Delta,\lambda,\iota,\zeta)$ is a $ccv$-graph if and only if all the following conditions are satisfied:
\begin{enumerate}
\item $\gcd(\lambda(x),\lambda(x^{-1})) = 1$ for all $x \in \D(\Delta)$;
\item $\zeta(x) \not \equiv \zeta(y)$ $(\mod \gcd(\iota(\beg x), \iota(\term x))$ for any two parallel darts $x,y \in \D(\Delta)$;
\item $\zeta(x) \not \equiv 0$ $(\mod \iota(\beg x))$ for all darts $x$ in a semi-edge; 
\item $\gcd(A\cup B) = 1$;
\item $\deg_{\lambda}(v) = 3$ for all $v \in \V(\Delta)$.
\end{enumerate}
\end{lemma}

In the lemma above, conditions (1)--(3) are there to guarantee the covering graph $\Cov(\Delta,\lambda,\iota,\zeta)$ is a simple graph, condition (4) that $\Cov(\Delta,\lambda,\iota,\zeta)$ is connected, and condition (5), that it is $3$-valent.

\subsection{Extendability to ccv-graphs}

Let $\Delta$ be a connected finite graph, and let $\lambda: \D(\Delta) \to \NN$ be an arbitrary function. We call the pair $(\Delta, \lambda)$ a {\em labelled graph}. Naturally, we obtain a labelled graph $(\Delta, \lambda)$ from every cyclic generalised voltage graph $(\Delta,\lambda,\iota,\zeta)$ by simply disregarding the functions $\iota$ and $\zeta$. In this case, the fact that $\lambda$ comes from a generalised voltage graph, which by definition satisfies equality (\ref{eq:ratio}), restricts $\lambda$ to some degree. We say that a labelled graph $(\Delta,\lambda)$ is {\em extendable} if there exist functions $\iota$ and $\zeta$ such that $(\Delta,\lambda,\iota,\zeta)$ is a cyclic generalised voltage graph. We will be particularly interested in those extendable labelled graphs that can be extended to a $ccv$-graph. 

A walk of length $n$ is a sequence of darts $W=(x_1,x_2,\ldots,x_n)$ such that $\beg x_{i+1} = \term x_i$ for all $i \in \{1,\ldots,n-1\}$. We say $W$ is {\em closed} if $\term x_n = \beg x_1$, and we say it is reduced if $x_{i+1} \neq x_i^{-1}$ for all $i \in \{1,\ldots,n-1\}$. A {\em path} is a reduced walk where $\beg x_i \neq \beg x_j$ for all $i \neq j$ and a {\em cycle} is a closed path. A cycle of length $n$ is also called an $n$-cycle. A tree is a connected graph without any cycles.

For a walk $W=(x_1,x_2,\ldots,x_n)$, we define the inverse of $W$ as the walk $W^{-1}=(x_n^{-1},x_{n-1}^{-1},\ldots,x_1^{-1})$. If $W_1=(x_1,x_2,\ldots,x_n)$ and $W_2=(y_1,y_2,\ldots,y_m)$ are two walks such that $\term x_n = \beg y_1$, then we define the concatenation of $W_1$ and $W_2$ as $W_1W_2 = (x_1,\ldots,x_n,y_1,\ldots,y_m)$.

Let $(\Delta,\lambda)$ be a labelled graph. We can extended the labelling $\lambda$ to a function $\lambda^*$ that assigns to each walk of $\Delta$ a rational number. For a walk $W=(x_1,x_2,\ldots,x_n)$ in $\Delta$ we let 
\begin{align}
\label{eq:rho}
 \lambda^*(W) := \prod\limits_{i=1}^n \frac{\lambda(x_i)}{\lambda(x_i^{-1})}.
\end{align}

We then have that
\begin{align}
\lambda^*(W^{-1}) = \lambda^*(W)^{-1} \text{ and } \lambda^*(W_1W_2)= \lambda^*(W_1)\lambda^*(W_2) 
\label{eq:rho2}
\end{align}
for any two walks $W_1$ and $W_2$ for which the concatenation is defined. The following is a useful characterization of labelled graphs that are extendable.

\begin{lemma}
{\rm \cite[Lemma 3.5]{cycliccovers}}
\label{lem:consistent}
A labelled graph $(\Delta,\lambda)$ is extendable if and only if $\lambda^*(C)=1$ for every closed walk $C$ of $\Delta$.
\end{lemma}

Naturally, not every extendable labelled graph $(\Delta,\lambda)$ can be extended to a $ccv$-graph. However, if $(\Delta,\lambda)$ does extend to a $ccv$-graph $(\Delta,\lambda,\iota,\zeta)$, we say $(\Delta,\lambda,\iota,\zeta)$ is a {\em $ccv$-extension} of $(\Delta,\lambda)$ and we call the covering graph $\Cov(\Delta,\lambda,\iota,\zeta)$ a {\em $ccv$-cover} of $(\Delta,\lambda)$.

\begin{lemma}
\label{lem:consistentccv}
\cite[Proposition 7.1]{cycliccovers}
A connected labelled graph  $(\Delta,\lambda)$ can be extended to a ccv-graph if and only if the following holds:
\begin{enumerate} 
\item $(\Delta,\lambda)$ is extendable;
\item $\lambda(x)=\lambda(x^{-1})$ implies $\lambda(x)=1$;
\item $\lambda(x) = \lambda(y) = 1$ for any two parallel darts $x$ and $y$;
\item $\lambda(x) = 1$ for every dart $x$ underlying a semi-edge;
\item $\deg_{\lambda}(v) = 3$ for all vertices $v \in \V(\Delta)$.
\end{enumerate}
\end{lemma}

Let $(\Delta, \lambda)$ be a labelled graph and let $x \in \D(\Delta)$ be such that $x \neq x^{-1}$ and $\lambda(x) \leq \lambda(x^{-1})$. We say $\{x,x^{-1}\}$ is an edge of type $[\lambda(x),\lambda(x^{-1})]$, or simply a $[\lambda(x),\lambda(x^{-1})]$-edge. From Lemma \ref{lem:consistentccv}, we see that if $(\Delta, \lambda)$ is extendable to a $ccv$-graph, then the edges of $\Delta$ are all of type $[1,1]$, $[1,2]$, $[1,3]$ or $[2,3]$
(see also \cite[Corollary 7.2]{cycliccovers}).

Now, suppose $\lambda(x) = 1$ for all $x \in \D(\Delta)$. That is, every edge of $(\Delta, \lambda)$ is a $[1,1]$-edge. By Lemmas \ref{lem:consistent} and \ref{lem:consistentccv}, $(\Delta, \lambda)$ admits a $ccv$-extension of $(\Delta,\lambda,\iota,\zeta)$ and a $ccv$-cover $\Gamma := \Cov(\Delta,\lambda,\iota,\zeta)$ . Since $\lambda(x) = 1$ for all $x \in \D(\Delta)$, it follows from formula (\ref{eq:ratio}) and the connectedness of $\Delta$ that $\iota(v) = \iota(u)$ for any two vertices $u$ and $v$ of $\Delta$. Then the canonical covering transformation $\rho$ of $\Gamma$ (mapping every vertex $v_i \in \V(\Gamma)$ to $v_{i+1}$) is a semiregular automorphism. Its vertex-orbits are precisely the vertex fibres of $\Gamma$. Therefore, if $\Delta$ has $k$ vertices, then $\rho$ has $k$ orbits on vertices and $\Gamma$ is a $k$-multicirculant graph.

\subsection{Simplified voltages}

\begin{figure}[h!]
\centering
\includegraphics[width=0.9\textwidth]{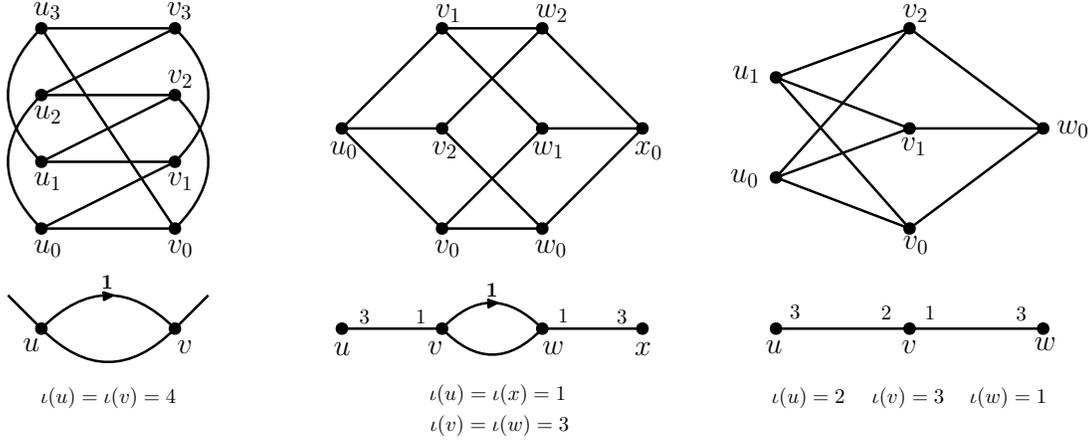}
\caption{A M\"{o}bius ladder, the cube graph and $K_{3,3}$ as covers of $ccv$-graphs (where the voltage is simplified).}
\label{fig:examples}
\end{figure}

As one would expect, for given a labelled graph $(\Delta,\lambda)$, there may exist different ccv-extension $(\Delta,\lambda,\iota,\zeta)$ and $(\Delta,\lambda,\iota,\zeta')$ such that $\Cov(\Delta,\lambda,\iota,\zeta) \cong \Cov(\Delta,\lambda,\iota,\zeta')$. It would be convenient to take, among all the possible $ccv$-extensions yielding isomorphic covers, one with a voltage assignment that is as ``nice'' as possible. As was proved  in \cite[Lemma 7.4]{cycliccovers}, for every $ccv$-extension $(\Delta,\lambda,\iota,\zeta)$ of a labelled graph $(\Delta,\lambda)$
there exists a voltage assignment $\zeta'$ that is {\em simplified} in the sense of Definition~\ref{def:simplzeta} below, such that
$\Cov(\Delta,\lambda,\iota,\zeta) \cong \Cov(\Delta,\lambda,\iota,\zeta')$.


\begin{definition}
\label{def:simplzeta}
Let $(\Delta,\lambda,\iota,\zeta)$ be a $ccv$-graph. The voltage $\zeta$ is a {\em simplified voltage} if for all $x \in \D(\Gamma)$ the following holds:
\begin{enumerate}
\item $\zeta$ is $\mathcal{T}$-normalised for a spanning tree $\mathcal{T}$ containing all $[i,j]$-edges with $i \neq j$;
\item $\zeta(x) < \gcd( \iota(\beg x),  \iota(\beg x^{-1}))$;
\item $\zeta(x) = \iota(\beg x)/2$ whenever $x$ underlies a semi-edge;
\item $0 < \zeta(x)$ and $\zeta(x) \neq \iota(\beg x)/2$ whenever $x$ underlies a loop.
\end{enumerate} 
\end{definition}
\noindent
Note that
if we assume a voltage $\zeta$ satisfies (1) in the definition above, then item (2) is equivalent to 
\begin{itemize}
\item[(2')] $\zeta(x) <  \iota(\beg x)$.
\end{itemize}
Indeed, clearly, (2) implies (2'). Suppose $\zeta$ satisfies (1) and (2'). If for some $x \in \D(\Gamma)$ we have $\zeta(x) = 0$ then $\zeta(x) <  \iota(\beg x)$. If $\zeta(x) \neq 0$, then by (1) we have $\iota(\beg x) =  \iota(\beg x^{-1})$ and thus $\gcd( \iota(\beg x),  \iota(\beg x^{-1})) = \iota(\beg x)$. Then $\zeta(x) < \iota(\beg x)=\gcd( \iota(\beg x),  \iota(\beg x^{-1}))$ and (2) holds.

\begin{remark}
One of the advantages of considering $ccv$-graphs with simplified voltage assignments is that the adjacency rules for the corresponding covering graphs become quite straightforward. Indeed, suppose $(\Delta,\lambda,\iota,\zeta)$ is a $ccv$-graph where $\zeta$ is simplified. Let $x \in \D(\Delta)$, let $u = \beg x$ and $v = \term x$, and let $e = \{x,x^{-1}\}$. Then for all $u_i \in \fib(u)$ we have:
\begin{enumerate}
\item If $u \neq v$ and $e$ is a $[1,1]$-edge, then $u_i \sim v_{i + \zeta(x)}$;
\item If $u \neq v$ and $e$ is a $[1,j]$-edge, $j \neq 1$, then $u_{i + k\cdot\iota(v)} \sim v_i$, with $0 \leq k < j$;
\item If $u = v$ and $e$ is a loop, then $u_i \sim u_{i+\zeta(x)}$ and $u_i \sim u_{i-\zeta(x)}$;
\item If $e$ is a semi-edge, then $u_i \sim u_{i+(\iota(u)/2)}$.
\end{enumerate}
\label{rem:adjacency}
\end{remark}


Let $(\Delta,\lambda,\iota,\zeta)$ be a $ccv$-graph. By \cite[Proposition 7.4]{cycliccovers} there exists a simplified voltage assignment $\zeta'$ for $\Delta$ such that $\Cov(\Delta, \lambda,\iota,\zeta) \cong \Cov(\Delta, \lambda,\iota,\zeta')$. This implies that every cubic graph is the cover of a $ccv$-graph with a simplified voltage assignment. For this reason, we will henceforth always assume that the voltage of a $ccv$-graph is simplified. The following theorem summarises the contents of this section.

\begin{theorem}
\label{the:coversum}
A graph $\Gamma$ is cubic and admits a cyclic group of automorphisms of order $n$ if and only if it is the cover of a cyclic generalised voltage graph $(\Delta,\lambda,\iota,\zeta)$ where $\lambda: \D(\Delta) \to \{1,2,3\}$ satisfies conditions (1)--(4) of Lemma \ref{lem:consistentccv}, $\zeta$ is simplified, $\gcd\{\zeta(x),\iota(v) \colon x \in \D(\Delta), v \in \V(\Delta)\}=1$ and $\lcm\{\lambda(x)\iota(\beg x) \colon x \in \D(\Delta)\}=n$.
\end{theorem}

\section{Vertex Transitive Covers}
\label{sec:VT}

We have shown in Section \ref{sec:covers} that if $\Gamma$ is a cubic graph admitting a cyclic group of automorphisms $G$, then $\Gamma$ is isomorphic to a $ccv$-cover of some labelled graph $(\Delta,\lambda)$, where the labelling $\lambda$ satisfies conditions (1)--(5) of Lemma \ref{lem:consistentccv}; that is, $(\Delta,\lambda)$ is extendable, $\lambda(x)=1$ for every dart $x$ underlying a loop, a semi-edge or a link that is parallel to another link, and $\deg_{\lambda}(x)=3$ for all darts $x \in \D(\Delta)$. If in addition we suppose that $\Gamma$ is vertex-transitive and that $G$ has order at least $\V(\Gamma)/3$, then further restrictions are set on the labelled graph $(\Delta,\lambda)$. It was shown in \cite{regorbs} that if $\Gamma$ is a cubic vertex-transitive graph and $G \leq \Aut(\Gamma)$ is cyclic, then the number of orbits of $G$ is bounded by a function of $k=|\V(\Gamma)|/|G|$. This, in turn, bounds the number of vertices of the quotient $\Gamma/G$. Furthermore, the ratio between the sizes of the largest and smallest orbits of $G$ is also bounded, which restricts the labelling $\lambda$. Let us be more precise.



Let $\Gamma$ be a cubic graph of order $n > 20$ admitting a cyclic subgroup of automorphisms $G$. Then $\Gamma$ is isomorphic to the cover of some $ccv$-graph $(\Delta,\lambda,\iota,\zeta)$ where $\lcm\{\lambda(x)\iota(\beg x)\colon x \in \D(\Delta) \}=|G|$. We can slightly abuse the language and identify $\Gamma$ with $\Cov(\Delta,\lambda,\iota,\zeta)$. Suppose $\Gamma$ is vertex-transitive and let $m = |G|$. Then by \cite[Theorem 4.7]{regorbs}, a $G$-orbit on vertices must have size $m/i$ for some $i \in \{1,2,3,4,6\}$ and the largest $G$-orbit has size precisely $m$. It follows that for any two vertices $\bar{u}$ and $\bar{v}$ of $\Gamma$, we have 
\begin{align*}
\frac{1}{6}|\bar{u}^G| \leq |\bar{v}^G| \leq 6|\bar{u}^G|.
\end{align*}
Since the orbits of $G$ on $\V(\Gamma)$ are identified with the fibres of $\V(\Delta)$ we see that 
\begin{align}
\frac{1}{6}\iota (u) \leq \iota(v) \leq 6\cdot\iota(u)
\label{eq:iotaineq}
\end{align}
for any two $u,v \in \V(\Delta)$ (recall that $\iota(u) = |\fib(u)|$ for all $u \in \V(\Delta)$). If in addition we suppose that $m \geq n/3$, then there must exist a vertex $\hat{u} \in \V(\Delta)$ such that

\begin{align}
3 \cdot \iota(\hat{u}) \geq \sum\limits_{v \in \V(\Delta)} \iota(v) = |\V(\Gamma)|.
\label{eq:iotathird}
\end{align}

Furthermore, \cite[Theorem 1.6]{regorbs} asserts that if $\Gamma$ has order $n>20$ and $\bar{u} \in \Gamma$ is such that $\bar{u}^G$ is of maximal order among all $G$-orbits, then $\bar{u}$ has a neighbour $\bar{v}$ such that $\bar{u}^G \neq \bar{v}^G$ but $|\bar{u}^G| = |\bar{v}^G|$. In particular, since the largest orbit of $G$ has size at least $n/3$, this implies that at least two thirds of the vertices of $\Gamma$ are contained in only two orbits ($\bar{u}^G$ and $\bar{v}^G$). It is an easy exercise to see that the remaining vertices can be divided in, at most, $3$ different orbits (of size $n/9$ each, or all three having different sizes: $n/6$, $n/9$ and $n/18$). Hence we obtain the following lemma.

\begin{lemma}
\label{lem:bound5}
Let $\Gamma$ be a cubic vertex transitive graph of order $n > 20$, and let $G \leq \Aut(\Gamma)$ be a cyclic group with an orbit of size $n/3$ or greater. Then $G$ has at most $5$ orbits on vertices.
\end{lemma}

Now, consider the function $\lambda^*$ defined with formula (\ref{eq:rho}). If $W=(x_1,x_2,\ldots,x_n)$ is a $uv$-walk in $\Delta$, then by a consecutive application of equality (\ref{eq:ratio}) to the darts of $W$ we have
\begin{align}
\iota(v) = \lambda^*(W)\iota(u).
\label{eq:iotawalk}
\end{align}
Since $u$ and $v$ are arbitrary vertices, we see that the index function $\iota$ is completely determined by the labelling $\lambda$ (as $\lambda^*$ depends only on $\lambda$) and the value of $\iota$ on a single vertex. The following proposition sums up the preceding paragraphs.

\begin{theorem}
\label{the:diagram}
Let $(\Delta,\lambda)$ be a labelled graph, let $\Gamma$ be a $ccv$-cover of $(\Delta,\lambda)$ of order $n>20$ and let $T$ be a spanning tree of $\Delta$. If $\Gamma$ is vertex-transitive and admits an automorphism of order $m \geq \frac{n}{3}$ then the following holds:
\begin{enumerate}
\item $(\Delta,\lambda)$ is extendable;
\item $\deg_{\lambda}(v) = 3$ for all vertices $v \in \V(\Delta)$;
\item $\lambda(x)=\lambda(x^{-1})$ implies $\lambda(x)=1$;
\item $\lambda(x) = \lambda(y) = 1$ for any two parallel darts $x$ and $y$;
\item $\lambda(x) = 1$ for every dart $x$ underlying a semi-edge;
\item $\Delta$ has at most $5$ vertices,
\end{enumerate}
moreover, there exists a vertex $\hat{u} \in \V(\Delta)$ such that: 
\begin{enumerate}
\item[(7)] $\hat{u}$ is incident to an edge of type $[1,1]$;
\item[(8)] $\frac{1}{6} \leq \lambda^*(W_v)$;
\item[(9)] $\sum\limits_{v \in \V(\Delta)\setminus \{\hat{u}\}} \lambda^*(W_v) \leq 2$;
\end{enumerate} 
for every $v \in \V(\Delta)\setminus\{\hat{u}\}$, where $W_v$ denotes the unique $\hat{u}v$-path in $T$.
\end{theorem} 

\begin{proof}
That items (1)--(5) hold follows at once from Lemma \ref{lem:consistentccv}. Item (6) holds by Lemma \ref{lem:bound5}. Now, let $u \in \V(\Gamma)$ be such that $|u^G|$ has maximum cardinality amongst all the vertex orbits of $G$. Then $|u^G| \geq \frac{n}{3}$. Let $\hat{u}=\pi(u)$ and let $v \in \V(\Delta)\setminus\{\hat{u}\}$ (recall that the natural projection $\pi$ maps every vertex $v_i \in \fib(v)$ to $v$). That $\hat{u}$ is incident to a $[1,1]$-edge follows from \cite[Theorem 4.7]{regorbs}, and thus (7) holds. Furthermore, by (\ref{eq:iotaineq}) we have $\frac{1}{6}\iota(\hat{u}) \leq \iota(v)$, but by (\ref{eq:iotawalk}) we can replace $\iota(v)$ by $\lambda^*(W_v)\iota(\hat{u})$. Thus (8) holds. To see that (9) holds, subtract $\iota(\hat{u})$ on both sides of inequality (\ref{eq:iotathird}) and replace $\iota(v)$ by $\lambda^*(W_v)\iota(\hat{u})$.
\end{proof}

\section{Artefacts}
\label{sec:obstructions}

A labelled subgraph of a labelled graph $(\Delta,\lambda)$ is a pair $(\Delta',\lambda')$ where $\Delta'$ is a subgraph of $\Delta$ and $\lambda'$ is the restriction $\lambda \mid_{\D(\Delta')}$. In this section we will define a set of `forbidden subgraph' for a labelled graph, that we call artefacts. An artefact in a labelled graph $(\Delta, \lambda)$ is a labelled subgraph of $(\Delta, \lambda)$ that guaranties the existence of a particular subgraph (containing a short cycle) in any  $ccv$-cover of $(\Delta, \lambda)$. 
We will show that a labelled graph containing certain  artefacts cannot admit a vertex-transitive $ccv$-cover of order large than $10$. First, we will define the notion of the signature of a graph.


Let $\Gamma$ be a cubic graph, let $x$ be a dart of $\Gamma$ and $c$ be a positive integer. Denote by $\epsilon_c(x)$ the number of $c$-cycles (cycles of length $c$) that pass through $x$. Let $v \in \V(\Gamma)$ and let $\{x_1,x_2,x_3\}$ be the set of darts beginning at $v$, ordered in such a way that $\epsilon_c(x_1) \leq \epsilon_c(x_2) \leq \epsilon_c(x_3)$. The triplet $(\epsilon_c(x_1),\epsilon_c(x_2),\epsilon_c(x_3))$ is then called the {\em $c$-signature} of $v$. Informally, the $c$-signature of $v$ tells us how the cycles of length $c$ passing through $v$ are distributed among the darts incident to $v$. If all vertices of $\Gamma$ have the same $c$-signature, we say that $\Gamma$ is {\em $c$-cycle-regular}, and we say the $c$-signature of $\Gamma$ is the $c$-signature of any of its vertices. 
Observe that if $\Gamma$ is vertex-transitive, then $\Gamma$ is $c$-cycle-regular for all $c \in \NN$. For cubic graphs of small girth $g$, the $g$-signature is sometimes enough to completely determine the graph.
The following lemma is a direct consequence of the results proved in \cite{signature} (or  
independently in \cite{g6alter}) and
\cite{girth6VT} .
\begin{lemma}
{\rm (See \cite[Theorem 1.5]{signature} and \cite[Theorem 1]{girth6VT}.)}
\label{lem:girthregular}
Let $\Gamma$ be a cubic girth-regular graph of girth $g \leq 6$. Then either the $g$-signature of $\Gamma$ is $(0,1,1)$ or one of the following occurs:
\begin{enumerate}
\item $g=3$ and $\Gamma \cong K_4$;
\item $g=4$ and one the following occurs
	\begin{enumerate}
	\item $\Gamma$ has signature $(1,2,2)$ and is isomorphic to a prism or a M\"{o}bius ladder;
	\item $\Gamma$ is isomorphic to $K_{3,3}$ or the cube graph $Q_3$;
	\end{enumerate} 
\item $g=5$ and $\Gamma$ is isomorphic to the Petersen graph or the dodecahedron $\GP(10,2)$.
\item $g=6$ and one of the following occurs:
    \begin{enumerate}
       \item  $\Gamma$ belongs to a finite list of exceptional graphs with at most $20$ vertices;
       \item $\Gamma$ has signature $(1,1,2)$, $(2,2,2)$ or $(3,4,5)$;
       \item $\Gamma$ has signature $(2,3,3)$ and is isomorphic either to
         \begin{itemize}
            \item a cyclic Haar graph
                   $\rH(3m;k,m)$ of order $6m$, $m>3$,
                        where $k=1$ if $m \equiv 0\> (\mod 3)$ and $k=3$ otherwise;
            \item a graph $\SDW(m,3)$ of order $6m$, $m>3$.
          \end{itemize}  
    \end{enumerate} 
\end{enumerate}
\end{lemma}

\begin{remark}
\label{rem:g6}
The Haar graph $\rH(3m;k,m)$ featuring in part (c) of the case $g=6$ of Lemma~\ref{lem:girthregular} was defined in \cite[Section 2.4]{girth6VT} as the bipartite Cayley graphs on the dihedral group $D_{3m} = \langle \rho, \tau \mid \rho^{3m}, \tau^2, (\rho\tau)^2\rangle$ with respect to the connection set $\{\tau, \tau\rho^k, \tau\rho^m\}$
and were denoted $\Delta_m$. It is, however, clear that such defined graph $\Delta_m$  and 
the cyclic Haar graph $\rH(3m;k,m)$ are isomorphic.
 As was proved there, their automorphism group has order $6m$ and thus acts regularly on the vertices.

 Similarly, the graph
$\SDW(m,3)$ was denoted in \cite{girth6VT} as $\Sigma_m$ and defined as the Cayley graph on the group $D_m\times \ZZ_3$
 with respect  to the connection set $S=\{(\rho\tau,0),(\tau,1), (\tau,2)\})$,
 where $D_m = \langle \rho, \tau \mid \rho^m, \tau^2, (\rho\tau)^2\rangle$ denotes the dihedral group of order $2m$. To prove that $\Sigma_m$ is indeed isomorphic to $\SDW(m,3)$ one can check that the mapping $(\rho^i\tau^j,k) \mapsto (i,k,1-j)$
 is indeed a graph isomorphism. As was proved in \cite[Proposition 5]{girth6VT},
  the automorphism group of $\SDW(m,3)$, $m>3$,
 is isomorphic to the group $S_3 \times D_m$ (where by $S_3$ we denote the symmetric group of order $6$).
\end{remark}

\begin{corollary}
\label{cor:smallgirthat}
If $\Gamma$ is a cubic arc-transitive graph of girth smaller than $6$, then $\Gamma$ is isomorphic to one of the following: $K_4$, $K_{3,3}$, the three-dimensional cube $Q_3$, the Petersen Graph or the dodecahedron $\rm{GP}(10,2)$.
\end{corollary}

\begin{lemma}\cite[Lemma 4.2]{ATgirth}
\label{lem:smallgirth2}
If $\Gamma$ is a cubic arc-transitive graph of girth $6$, then either $\Gamma$ has $6$-signature $(2,2,2)$ or it has order $n \leq 20$.
\end{lemma}

\begin{lemma}
\label{lem:13edge}
Let $(\Delta, \lambda,\iota,\zeta)$ be a $ccv$-graph and suppose $\Gamma:= \Cov(\Delta, \lambda,\iota,\zeta)$  is vertex-transitive. If for some $x \in \D(\Delta)$  we have $\lambda(x)=3$, then $\Gamma$ is arc-transitive.
\end{lemma}

\begin{proof}
Let $u = \beg x$ and $u_0 \in \fib(u)$. Let $a$ and $b$ be two arbitrary darts of $\Gamma$. Since $\Gamma$ is vertex-transitive, there exist automorphisms $\phi,\psi \in \Aut(\Gamma)$ such that $(\beg a)^{\phi}=u_0=(\beg b)^{\psi}$. In particular, both $a^{\phi}$ and $b^{\psi}$ begin at $u_0$. Furthermore, since $\lambda(x)=3$, all three darts beginning at $u_0$ belong to the same orbit of the cyclic group of automorphisms of $\Gamma$ preserving the fibres, and thus there exists $\gamma \in \Aut(\Gamma)$ such that $a^{\phi\gamma} = b^{\psi}$. Then $a^{\phi\gamma\psi^{-1}}=b$. We conclude $\Gamma$ is arc-transitive.
\end{proof}

\begin{figure}[h!]
\centering
\includegraphics[width=0.8\textwidth]{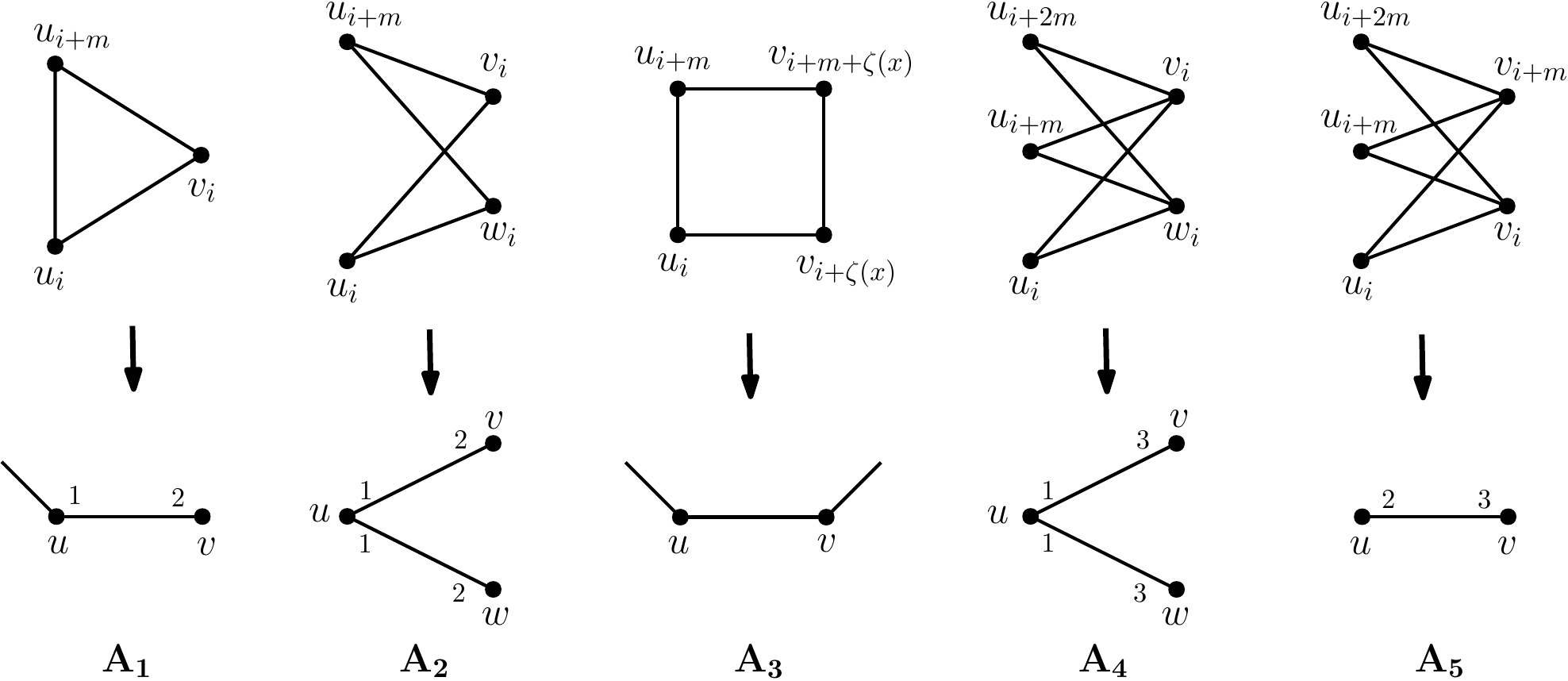}
\caption{The five artefacts $A_i$ with $i \in \{1,2,3,4,5\}$ (bottom row). Above each, a small subgraph of the cover of any $ccv$-graph containing $A_i$.}
\label{fig:A15}
\end{figure}

Let $A_1$, $A_2$, $A_3$, $A_4$ and $A_5$ be the $5$ labelled graphs depicted in the bottom row of Figure \ref{fig:A15}.  

\begin{lemma}
\label{lem:artcycles1}
Let $(\Delta,\lambda)$ be a labelled graph and let $\Gamma$ be a $ccv$-cover of $(\Delta,\lambda)$. Suppose $(\Delta,\lambda)$ contains an artefact $A_j$ for some $j \in \{1,2,3,4,5\}$. Then: 
\begin{enumerate}
\item if $j=1$, $\Gamma$ contains a $3$-cycle;
\item if $j \in \{2,3\}$, $\Gamma$ contains a $4$-cycle;
\item if $j \in \{4,5\}$, $\Gamma$ contains a copy of $K_{3,2}$.
\end{enumerate}
\end{lemma}
\begin{proof}
The proof consists in repeatedly applying Remark \ref{rem:adjacency} for each of the five possible cases. For instance, suppose $\Gamma \cong \Cov(\Delta,\lambda,\iota,\zeta)$ where $(\Delta,\lambda,\iota,\zeta)$ is a $ccv$-extension of $(\Delta,\lambda)$ and that $\Gamma$ contains a subgraph isomorphic to $A_1$. Assume the notation of Figure \ref{fig:A15}. Then $\iota(u)=2\cdot\iota(v)=2m$ for some $m \in \NN$ and since $uv$ is a $[1,2]$-edge, it follows from Remark \ref{rem:adjacency} that each $v_i \in \fib(v)$ is adjacent to $u_i$ and $u_{i+m}$. Furthermore, since $u$ is incident to a semi-edge, $u_i$ is incident to $u_{i+m}$ (again, by Remark \ref{rem:adjacency}). Then, $(u_i,u_{i+m},v_i)$ is a $3$-cycle of $\Gamma$ for all $i \in \{0,\ldots,m-1\}$. A similar argument shows that the lemma holds for the four remaining cases (see Figure \ref{fig:A15}).
\end{proof}

\begin{corollary}
\label{cor:ATart}
Let $(\Delta,\lambda)$ be labelled graph containing an artefact $A_i$ with $i \in \{1,2,3,4,5\}$ and a dart $x$ such that $\lambda(x)=3$. If $\Gamma$ is a vertex-transitive $ccv$-cover of $(\Gamma,\lambda)$, then $\Gamma$ is isomorphic to $K_4$, $K_{3,3}$ or $Q_3$.  
\end{corollary}

\begin{proof}
By Lemma \ref{lem:13edge} $\Gamma$ is arc-transitive and since it contains an artefact $A_i$ with $i \in \{1,2,3,4,5\}$ it has girth $3$ or $4$. It follows from Corollary \ref{cor:smallgirthat} that $\Gamma$ is isomorphic to $K_4$, $K_{3,3}$ or $Q_3$.
\end{proof}

\begin{figure}[h!]
\centering
\includegraphics[width=0.6\textwidth]{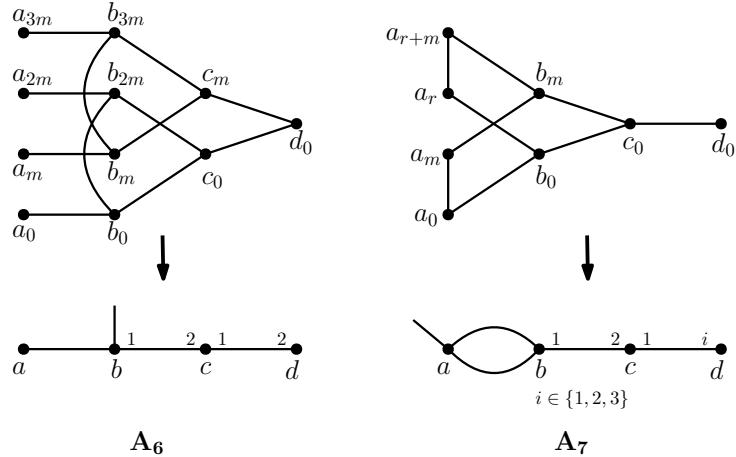}
\caption{The two artefacts $A_6$ and $A_7$ (bottom row). Above each, a small subgraph of the cover of any $ccv$-graph containing it.}
\label{fig:A69}
\end{figure}

Let $A_6$ and $A_7$ be the labelled graphs depicted in the bottom row of Figure \ref{fig:A69}.

\begin{lemma}
If $(\Delta,\lambda)$ is a labelled graph containing $A_6$ or $A_7$, then no $ccv$-cover of $(\Delta,\lambda)$ is vertex-transitive.
\label{lem:A9}
\end{lemma}

\begin{proof}
Let $(\Delta,\lambda,\iota,\zeta)$ be a $ccv$-extension of $(\Delta,\lambda)$ such that $\Gamma\cong \Cov(\Delta,\lambda,\iota,\zeta)$. Suppose, for a contradiction, that $\Gamma$ is vertex-transitive and $(\Delta,\lambda)$ contains a copy of $A_j$ for some $j \in \{6,7\}$.

Now, suppose $j=7$. Let $a$, $b$, $c$ and $d$ be the vertices of $A_7$, as they are labelled in Figure \ref{fig:A69}. There are two edges connecting $b$ to $a$. Without loss of generality, we may assume the darts on one of these edges have trivial voltage, as necessarily one of these edges must lie on a spanning tree of $\Delta$ and (we can assume) the voltage assignment $\zeta$ is simplified. As for the other edge, let $r \neq 0$ be the voltage of the dart underlying it and beginning at $b$ (and thus, its inverse, beginning at $a$, has voltage $-r$). Let $\iota(c)=m$ so that $\iota(b)=\iota(a)=2m$. Observe that $(c_0,b_0,a_0,a_m,b_m)$ and $(c_0,b_0,a_{r},a_{m+r},b_m)$ are $5$-cycles of $\Gamma$ (note that this is true even if $r=m$). Then every dart beginning at $b_0$ lies on a $5$-cycle. Since $dc$ is a $[1,i]$-edge of $(\Delta,\lambda)$, we see that $d_0$ is adjacent to $c_0$ in $\Gamma$. Furthermore, since $\Gamma$ is vertex-transitive, the dart beginning at $d_0$ and ending at $c_0$ lies on a $5$-cycle $C:=(d_0,c_0,u,v,w)$, for some $u,v,w \in \V(\Gamma)$. Clearly, $u \in \{b_0,b_m\}$ and $v \in \{a_0,a_r,a_m,a_{m+r}\}$, since both $\fib(c)$ and $\fib(b)$ are independent sets. Then $w \in \fib(a)\cup\fib(b)$ which leads us to a contradiction, since no vertex in $\fib(a)\cup\fib(b)$ is adjacent to a vertex in $\fib(d)$. Therefore, there is no $5$-cycle tracing the edge $d_0c_0$ and thus $\Gamma$ is not vertex-transitive.

Finally, suppose $j=6$ and assume the notation in Figure \ref{fig:A69}. Then $\iota(d) = m$ for some $m \in \NN$. By Lemma \ref{lem:artcycles1}, $\Gamma$ contains a $3$-cycle and thus the vertex $d_0 \in \fib(d)$ must lie on a $3$-cycle $C$. Since $d_0$ has two neighbours in $\fib(c)$, one vertex of $C$ must be in $\fib(c)$. Without loss of generality, let $c_0$ be that vertex. Then, the third vertex in $C$ must be a common neighbour of $d_0$ and $c_0$, but the other two neighbours of $c_0$ are $b_0$ and $b_{2m}$, none of which si adjacent to $d_0$. Therefore, $d_0$ does not lie on a $3$-cycle, and $\Gamma$ is not vertex-transitive.
\end{proof}

The following Theorem, which is a consequence of Corollary \ref{cor:ATart} and Lemma \ref{lem:A9}, summarises the contents of this section.

\begin{theorem}
Let $(\Delta,\lambda)$ be a labelled graph and let $\Gamma$ be a ccv-cover of $(\Delta,\lambda)$ with more than $20$ vertices. If $\Gamma$ is vertex transitive, then one of the following hold:
\begin{enumerate}
\item $(\Delta,\lambda)$ does not contain an artefact $A_i$ with $i \in \{1,2,3,4,5\}$ and a dart $x$ such that $\lambda(x)=3$;
\item $(\Delta,\lambda)$ does not contain an artefact $A_i$ with $i \in \{6,7\}$.
\end{enumerate}
\label{the:artefacts}
\end{theorem}

\section{The set $\mQ$}
\label{sec:deltas}

Recall that $\mG$ is the set of all vertex-transitive cubic graphs $\Gamma$ of order $n>20$ admitting an automorphism $g$ of order $n/3$ or greater, and that $\mQ$ is the set of all labelled quotients $\Gamma / \langle g \rangle$. As a step towards proving Theorem \ref{the:main}, we must determine the set $\mQ$. Then, we can reconstruct $\mG$ by considering all vertex-transitive $ccv$-covers of elements of $\mQ$. Observe that $\mQ$ is a subset of the set $\mQ^*$ of all labelled graphs satisfying the conditions stated in Theorems \ref{the:diagram} and \ref{the:artefacts}. These conditions are restrictive enough to allow us to quickly compute $\mQ^*$ by means of a brute-force algorithm. As it transpires, $\mQ^*$ consist of $20$ labelled graphs. To determine $\mQ$ it suffices to determine which of these graphs admit a vertex-transitive $ccv$-cover with more than $20$ vertices. 

\begin{figure}[h!]
\centering
\includegraphics[width=0.6\textwidth]{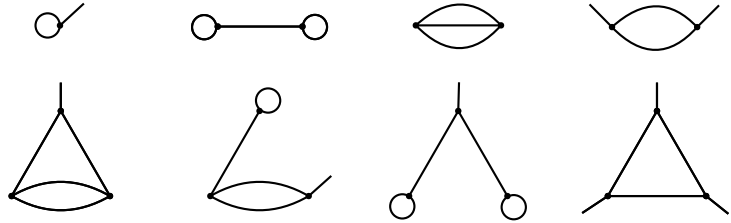}
\caption{Elements of $\mQ^*$ with at most $3$ vertices.}
\label{fig:kcirc}
\end{figure}

The eight elements of $\mQ^*$ having less than four vertices, shown in Figure \ref{fig:kcirc}, correspond to the eight possible quotients of a cubic graph by a $(k,\ell)$-semiregular automorphism with $k \in \{1,2,3\}$ (graphs with such an automorphism are also called $k$-multicirculants).

 Each of these eight graphs admit at least one vertex-transitive $ccv$-cover with more than $20$ vertices (see \cite[Theorem 1.1]{cycliccovers} or \cite{bic} and \cite{tricirc} for details). Therefore, the eight graphs of Figure \ref{fig:kcirc} are elements of $\mQ$.

\begin{figure}[h!]
\centering
\includegraphics[width=1\textwidth]{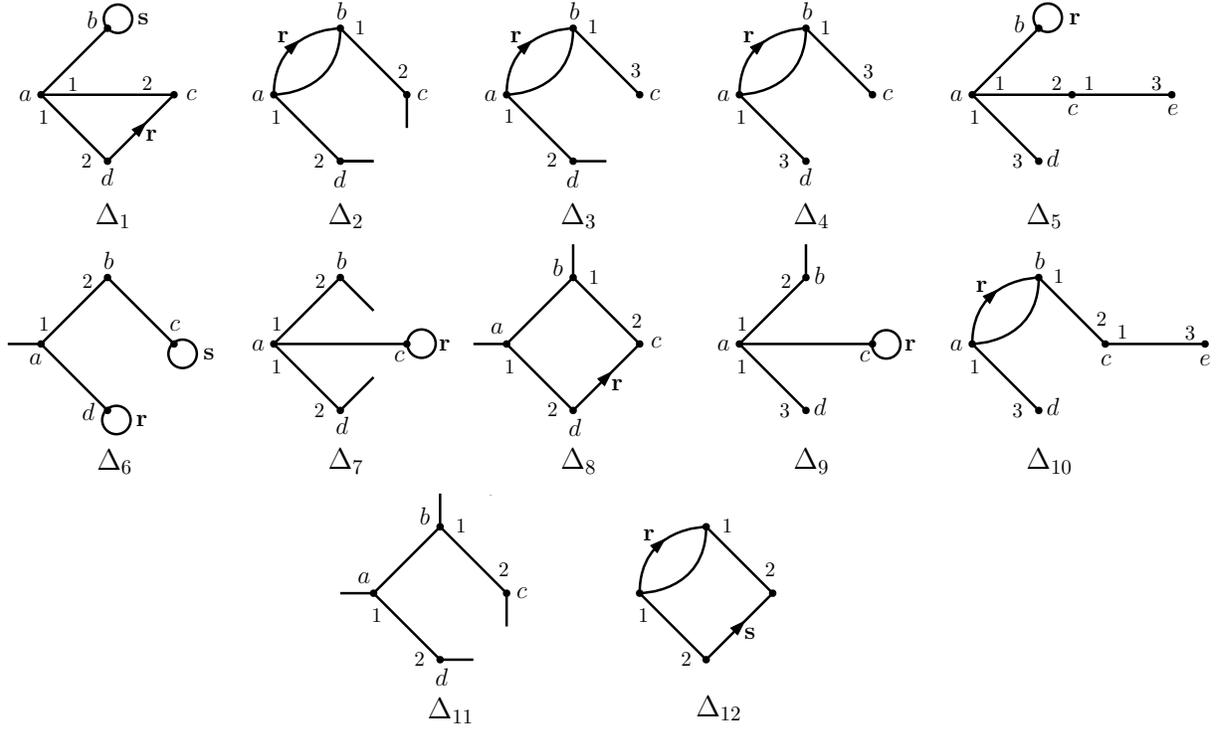}
\caption{Elements of $\mQ^*$ with more than $3$ vertices.}
\label{fig:exceptions}
\end{figure}

The remaining $12$ labelled graphs of $\mQ^*$ are shown in Figure \ref{fig:exceptions}. We will show that, with the exception of  $\Delta_{12}$, none of these graphs admit a vertex-transitive $ccv$-cover of order larger than $20$. The graph $\Delta_{12}$ will be studied in detail in Section \ref{sec:exception}.

Now consider a labelled graph $\Delta_i$ from Figure \ref{fig:exceptions} and suppose $(\Delta,\lambda,\iota,\zeta)$ is a $ccv$-extension of $\Delta_i$. We may assume that $\zeta$ is a simplified voltage assignment and agrees with Figure \ref{fig:exceptions}, where we adopt the following notation convention. For a symbol $\alpha \in \{r,s\}$, an edge with an arrow oriented from, say, $u$ to $v$, with the letter $\alpha$ next to it, indicates that the dart $x$ underlying this edge and beginning at $u$ has voltage $\alpha$, for some $0 \leq \alpha < \iota(v)$. A loop with the letter $\alpha$ next to it, indicates that one of the underlying darts, say $x$, has voltage $\alpha$ for some integer $0 < \alpha < \iota(\beg x)$. For a semi-edge $x$, $\zeta(x) = \iota(\beg x)/2$. All other darts belong to a spanning $\mathcal{T}$ and have trivial voltage. The vertices of each $\Delta_i$ are named in Figure \ref{fig:exceptions}, but we refrain from naming the darts in the figure so as not to overburden it. Since parallel darts in a $ccv$-graph need to have distinct voltages, every dart in $\Delta_i$ is completely determined by its endpoints along with its voltage. Hence, we will denote a dart $x$ of $\Delta_i$ beginning at $u$ and ending at $v$ by $(uv)_{\zeta(x)}$; its inverse is then $(vu)_{-\zeta(x)}$. As every cover of a $ccv$-graph is a simple graph, the darts in the fibre of $(uv)_{\zeta(x)}$ are denoted by $u_iv_{i+\zeta}$ for $i \in \{0,\ldots,\iota(u)-1\}$ (like arc or directed edges are usually denoted). 

Let $(\Delta,\lambda)$ be a labelled graph, $\Gamma$ be a $ccv$-cover of $(\Delta,\lambda)$ and let $\pi: \Gamma \to \Delta$ be the corresponding projection. If $W$ is a $uv$-walk in $\Delta$, then a {\em lift} of $W$ based at a vertex $u_i \in \fib(u)$ is a walk $\overline{W}=(x_1,x_2,\dots,x_n)$ beginning at $u_i$ such that the projection $\pi(\overline{W}):=(\pi(x_1),\pi(x_2),\ldots,\pi(x_n))$ is equal to $W$. We denote by $\mathcal{L}(W)$ the set of all lifts of $W$ based at $u_0$.

We say $W$ is $\lambda$-reduced if $x_{i+1} \neq x_i^{-1}$ whenever $\lambda(x_i^{-1}) = 1$, and $x_n \neq x_1^{-1}$ whenever $\lambda(x_1) = 1$. Clearly, every reduced walk is $\lambda$-reduced. 


Let $(\Delta,\lambda,\iota,\zeta)$ be a $ccv$-graph and set $n = \lcm \{\lambda(x)\iota(\beg x) \colon x \in \D(\Delta)\}$. Let $W=(x_1,x_2,\ldots,x_k)$ be a walk in $\Delta$ and let $d = \gcd\{\iota(\beg x_i)\colon x_i \in W\}$. We define the {\em endset} of $W$ as
$$ \textrm{end}(W) = \sum\limits_{i=0}^k \zeta(x_i) + \langle d \rangle,$$
where $\langle d \rangle$ denotes the subgroup of $\ZZ$ generated by $d$, and where the addition is computed modulo $\iota(\term x_k)$.

\begin{lemma}\cite[Lemma 30]{MPgc}
\label{lem:walks}
Let $(\Delta, \lambda,\iota,\zeta)$ be a $ccv$-graph and $\Gamma= \Cov(\Delta, \lambda,\iota,\zeta)$. Let $W$ be a $uv$-walk for some $u,v \in \V(\Delta)$. If $\overline{W} \in \mathcal{L}(W)$, then the final vertex of $\overline{W}$ is $v_{j}$ for some $j \in \term(W)$. Conversely, for every $j \in \term(W)$ there exists a lift of $W$ beginning at $u_0$ and ending at $v_j$.
\end{lemma}  

\begin{lemma}
\label{cor:walks}
Let $(\Delta, \lambda,\iota,\zeta)$ be a $ccv$-graph and $\Gamma= \Cov(\Delta, \lambda,\iota,\zeta)$. If $C$ is a cycle in $\Gamma$, then $\pi(C)$ is a $\lambda$-reduced closed walk in $\Delta$ and $0 \in \term(\pi(C))$.
\end{lemma}

\begin{proof}
Let $C$ be a cycle in $\Gamma$. Clearly, $C$ is a $u_au_a$-walk for some vertex $u_a \in \fib(u)$ and some $u \in \V(\Delta)$. 
Let $x$ be a dart visited by $\pi(C)$ and suppose that $\pi(C)$ traces $x^{-1}$ immediately after $x$. We will show that  $\lambda(x^{-1}) \neq 1$. Observe that since $\pi(C)$ traces $x$ and $x^{-1}$ one after the other, there must exist a dart $x_i \in \fib(x)$ and a dart $(x^{-1})_j \in \fib(x^{-1})$ such that $C$ traces $x_i$ and $(x^{-1})_j$ consecutively. Since $C$ is a cycle, it is a reduced walk by definition, and thus $(x^{-1})_j \neq (x_i)^{-1}$. However both $(x^{-1})_j$ and $(x_i)^{-1}$ belong to $\fib(x^{-1})$. Moreover, since $x_i$ and $(x^{-1})_j$ are two consecutive darts of a walk, we have  $\beg (x^{-1})_j = \term x_i = \beg (x_i)^{-1}$. That is, there are two distinct darts in $\fib(x^{-1})$ beginning at the same vertex. This implies that $\lambda(x^{-1}) \geq 2$. 

Now suppose $x$ and $x^{-1}$ are the first and the last darts traced by $\pi(C)$. Let $x_i \in \fib(x)$ and $(x^{-1})_j \in \fib(x^{-1})$ be the first and last darts traced by $C$, respectively. By an argument analogous to the one used in the previous case, we have that $(x^{-1})_j \neq x_i^{-1}$ but $\beg(x^{-1})_j = \beg x_i^{-1}$, and thus $\lambda(x)\geq 2$. This shows that $\pi(C)$ is $\lambda$-reduced.

To show that $0 \in \term(\pi(C))$, recall that $\Gamma$ admits an automorphism $\rho$ that maps every dart $x_i$ to $x_{i+1}$. Then $\rho^{-a}$ maps the vertex $u_a$ to $u_0$, and thus $\rho^{-a}(C)$ is a reduced closed walk beginning and ending at $u_0$. Moreover, $\rho^{-a}(C) \in \mathcal{L}(\pi(C))$. Since the final vertex of $\rho^{-a}(C)$ is $u_0$, it follows from Lemma \ref{lem:walks} that $0 \in \term(\pi(W))$.
\end{proof}

Throughout the rest of the section we will assume that $(\Delta,\lambda)=\Delta_i$ for some $i \in \{1,\ldots,11\}$, that $\Gamma$ is the cover of a $ccv$-extension $(\Delta,\lambda,\iota,\zeta)$ of $\Delta_i$, and that $\pi:\Gamma \to \Delta$ is the covering projection. Note that $\Gamma$ is completely determined by the values of the voltages $r$ and $s$, and $m := \gcd\{\iota(u) \colon u \in \V(\Delta)\}$.
Indeed, $\Gamma$ is uniquely determined by the quadruple $(\Delta,\lambda,\iota,\zeta)$. 
Recall that the index function $\iota$ is determined by its value on a single vertex along with the labelling $\lambda$, which is given. Let $u \in \V(\Delta)$ be any vertex and for each $v \in \V(\Delta) \setminus \{u\}$ chose (arbitrarily) a $uv$-walk $W_v$. By equality (\ref{eq:iotawalk}), we have $\iota(v)=\lambda^*(W_v)\iota(u)$ for all $v \in \V(\Delta) \setminus \{u\}$. Let $c$ be the smallest positive integer such that $c \cdot \lambda^*(W_v)$ is an integer for all $v \in \V(\Delta) \setminus \{u\}$. Then $\iota(u)=c\cdot m$. Note that $c$ depends only on $\lambda$ and our choice of $u$. Then, $\iota$ is completely determined by $\lambda$ and $m$. Finally, since we can assume $\zeta$ to be simplified, we know that every dart $x$ underlying a semi-edge has voltage $\iota(\beg x) / 2$, and any other dart not labelled $r$ or $s$ has voltage $0$. The values of $r$ and $s$, along with the function $\iota$ thus completely determine $\zeta$.

We are now ready to analyse the labelled graphs $\Delta_i$. The technique employed in the following pages relies mainly in finding a closed walk $W$ of length $n$ in $\Delta_i$ such that $\mathcal{L}(W)$ contains a cycle of length $n$, regardless of the specific values of the voltages $r$ and $s$. Such a walk can often be found by finding an artefact $A_j$ in $\Delta_i$. If we suppose that $\Gamma$ is vertex-transitive, then for every vertex $v_i$ of $\Gamma$, at least one dart incident to $v_i$ must lie on an $n$-cycle. Then by Corollary \ref{cor:walks}, this will imply that for every vertex $v \in \V(\Delta)$ a specific dart incident to $v$ lies on a closed walk $W'$ of length $n$ such that $0 \in \term(W')$. Since every element of $\term(W')$ can be seen as a linear combination of $m$, $r$ and $s$,  $0 \in \term(W')$ implies a relation between $m$, $r$ and $s$, which along with the fact that $\gcd(m,r,s)=1$ (see item (5) of Lemma \ref{lem:ccv}), is often enough to completely determine their values (up to a few options). This, in turn, determines the graph $\Gamma$.

\begin{lemma}
If $\Gamma$ is a vertex-transitive $ccv$-cover of $\Delta_{1}$, then $\Gamma$ is a bicirculant graph of order $12$.
\end{lemma}

\begin{proof}
Let $(\Delta,\lambda,\iota,\zeta)$ be a $ccv$-extension of $\Delta_1$ such that $\Gamma = \Cov(\Delta,\lambda,\iota,\zeta)$. Let $\iota(c)=m$ and note that $\iota(d)=m$ as $c$ and $d$ are connected through a $[1,1]$-edge. Similarly $\iota(a)=\iota(b)=2m$. Let $W=((ac)_0,(ca)_0,(ad)_0,(da)_0)$ and note that every reduced walk in $\mathcal{L}(W)$ is a cycle of length $4$ (see Lemma \ref{lem:artcycles1}). Then, every dart in the fibre of $(ac)_0$ or $(ad)_0$ lies on a $4$-cycle. Since $\Gamma$ is vertex-transitive, every dart in the fibre of $(bb)_s$ (or of $(bb)_{-s}$) must lie on a $4$-cycle. Suppose $C$ is a $4$-cycle (in $\Gamma$) through $b_0b_s$. Then $\pi(C)$ is a $\lambda$-reduced walk of length $4$ through $(bb)_s$. Clearly, $\pi(C) = ((bb)_s,(bb)_s,(bb)_s,(bb)_s)$ and $\term(\pi(C)) = \{4s\}$. By Lemma \ref{cor:walks}, $0 \in\term(\pi(C))$ and thus $4s \equiv 0$ ($\mod 2m$). Since $0 < s < m$, we have $2s = m$. This shows that $(b_0,b_s,b_{2s},a_{2s},c_{2s},a_{2s+m})$ is a $6$-cycle in $\Gamma$ since $a_{2s+m} = a_0$ and $a_0 \sim b_0$. Then the $6$-signature of $\Gamma$ is $(\epsilon_1,\epsilon_2,\epsilon_3)$ where $\epsilon_i > 0$. That is, every dart of $\Gamma$ lies in at least one $6$-cycle. In particular, there is a $6$-cycle $C'$ through $c_0d_{-r}$. Then $\pi(C')$ is a $\lambda$-reduced walk of length $6$ through $(cd)_{-r}$. By inspecting Figure \ref{fig:exceptions}, we see that $\pi(C')$ must be one of the following:
\begin{align*}
W_1 &= ((cd)_{-r},(da)_0,(ac)_0,(cd)_{-r},(da)_0,(ac)_0),\\
W_2 &= ((cd)_{-r},(da)_0,(ab)_0,(bb)_s,(ba)_0,(ac)_0),\\
W_3 &= ((cd)_{-r},(da)_0,(ab)_0,(bb)_{-s},(ba)_0,(ac)_0).
\end{align*}
Now, $\term(W_1)=\{-2r\}$, $\term(W_2)=\{s-r\}$ and $\term(W_3)=\{-s-r\}$. Since 
$0 \in \term(W_i)$ for some $i \in \{1,2,3\}$, we see that one of the following holds modulo $m$,
\begin{align*}
-2r \equiv 0,\\
s-r \equiv 0,\\
-s-r \equiv 0.
\end{align*}
Since $0 \leq r <m$ and $2s = m$, we see that either $r=0$ or $s=r$. However, $\gcd(m,r,s)=1$ by Lemma~\ref{lem:ccv}. Therefore, $r = s = 1$ and $m=2$. That is, the functions $\zeta$ and $\iota$ are completely determined and so is $\Gamma$. It can be verified that $\Gamma$ is a bicirculant isomorphic to the Franklin graph (see page 244 of \cite{bondy} for definition and properties).
\end{proof}

\begin{figure}[h!]
\centering
\includegraphics[width=0.7\textwidth]{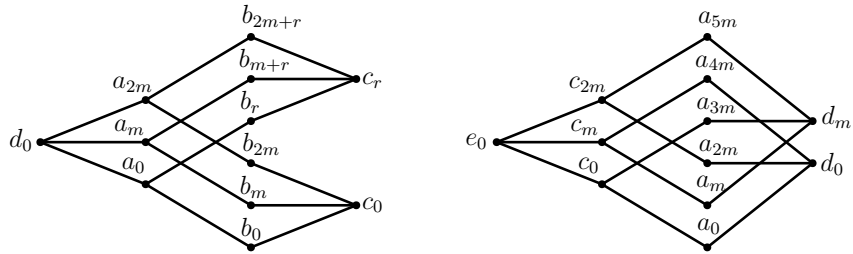}
\caption{For each of the two cases in the proof of Lemma \ref{lem:delta45}, a subgraph of $\Gamma$ containing $4$ distinct $6$-cycles through the dart $d_0a_0$ (left) and through the dart $e_0c_0$ (right), respectively.}
\label{fig:delta4delta5}
\end{figure}

\begin{lemma}
If $\Gamma$ is a $ccv$-cover of $\Delta_{i}$, $i \in \{2,3,11\}$, then $\Gamma$ is not vertex-transitive.
\end{lemma}

\begin{proof}
Let $(\Delta,\lambda,\iota,\zeta)$ be a $ccv$-extension of $\Delta_j$ and suppose $\Gamma := \Cov(\Delta,\lambda,\iota,\zeta)$ is vertex-transitive. 

First, suppose $j = 2$. Let $\iota(d) = 2m$ and for $k \in \{0,1\}$ let $W_k=( (da)_0,(ab)_{kr},(bc)_0,(cb)_0,(ba)_{-kr},(ad)_{0})$. Observe that every edge incident to $a_0$ lies on a $6$-cycle belonging to $\mathcal{L}(W_0) \cup \mathcal{L}(W_1)$. Then $d_0d_m$ must lie on a $6$-cycle $C$ of $\Gamma$ and $\pi(C)$ is a $\lambda$-reduced closed walk of length $6$. It is straightforward to see that no $\lambda$-reduced closed walk of length $6$ in $\Delta_2$ traces the dart $(dd)_m$, a contradiction. Therefore $\Gamma$ is not vertex-transitive.

Now, suppose $j = 3$. Let $\iota(d) = 2m$ and $W=( (da)_0,(ab)_r,(ba)_0,(ad)_0,(da)_0,(ab)_0,(ba)_{-r},(ad)_0)$. Observe that every edge incident to $a_0$ lies on an $8$-cycle of $\mathcal{L}(W)$, which implies the existence of an $8$-cycle $C$ through $d_0d_m$, since $\Gamma$ is vertex-transitive. Then $\pi(C)$ is a $\lambda$-reduced closed walk of length $8$ through $(dd)_m$. Once more, one can verify that no such walk exists in $\Delta_3$. We conclude $\Gamma$ is not vertex-transitive.

Finally, suppose suppose $j = 11$. Let $\iota(d) = m$ and note that since the voltage assignment is simplified, the darts incident to $a$ or $b$ have voltage $m$ and those incident to $c$ or $d$ have voltage $m/2$. Since $\Gamma$ is connected, $\gcd(m/2,m)=1$, which implies that $m=2$. Note that the functions $\zeta$ and $\iota$ are thus completely determined, and so is the graph $\Gamma$. One can simply verify that $\Gamma$ is a non-vertex-transitive graph of order $12$.
\end{proof}

\begin{lemma}
If $\Gamma$ is a vertex-transitive $ccv$-cover of $\Delta_j$, with $j \in \{4,5\}$, then $\Gamma$ has less than $20$ vertices.
\label{lem:delta45}
\end{lemma}

\begin{proof}
Suppose $\Gamma=\Cov(\Delta,\lambda,\iota,\zeta)$ is vertex-transitive where $(\Delta,\lambda,\iota,\zeta)$ is a $ccv$-extension of $\Delta_{4}$ and let $\iota(a)=m$. Since $\Delta_4$ has a $[1,3]$-edge, $\Gamma$ must be arc-transitive by Lemma \ref{lem:13edge}. Now, if $r \equiv 0$ ($\mod m$), then $(d_0,a_0,c_m,a_m)$ is a $4$-cycle of $\Gamma$ and by Corollary \ref{cor:smallgirthat}, $\Gamma$ has less than $20$ vertices. Suppose that $r \not \equiv 0$ ($\mod m$). For $i \in \{0,1\}$, let 
$$W_i=((da)_0,(ab)_{ir},(bc)_0,(cb)_0,(ba)_{-ir},(ad)_0).$$ 
Observe that every dart in the fiber of $(da)_0$ lies on $4$ distinct $6$-cycle in $\mathcal{L}(W_0) \cup \mathcal{L}(W_1)$ (see Figure \ref{fig:delta4delta5}, left). Then by Lemma \ref{lem:smallgirth2} $\Gamma$ has less than $20$ vertices.

Now, suppose $\Gamma:=\Cov(\Delta,\lambda,\iota,\zeta)$ is vertex-transitive where $(\Delta,\lambda,\iota,\zeta)$ is a $ccv$-extension of $\Delta_{5}$ and consider the walk 
$$W=((ec)_0,(ca)_0,(ad)_0,(da)_0,(ac)_0,(ce)_0).$$ 
Observe that every dart in  $e_0c_0$ lies on $4$ distinct $6$-cycle in $\mathcal{L}(W)$ (see Figure \ref{fig:delta4delta5}, right) and by Lemma \ref{lem:smallgirth2}, $\Gamma$ has less than $20$ vertices.
\end{proof}

\begin{lemma}
If $\Gamma$ is a vertex-transitive $ccv$-cover of $\Delta_{6}$, then $\Gamma$ is a tricirculant of order $18$.
\end{lemma}

\begin{proof}
Suppose $\Gamma:=\Cov(\Delta,\lambda,\iota,\zeta)$ is vertex-transitive where $(\Delta,\lambda,\iota,\zeta)$ is a $ccv$-extension of $\Delta_{6}$. Then, for some $m \in \ZZ$ we have $\iota(c)=\iota(d)=m$ and $\iota(a)=\iota(b)=2m$. Let $W=((aa)_m,(ad)_0,(da)_0)$ and observe that every reduced walk in $\mathcal{L}(W)$ is a $3$-cycle (see Lemma \ref{lem:artcycles1}). Then every dart in the fibre of $(aa)_m$ or $(ad)_0$ lies on a $3$-cycle.
Since $\Gamma$ is vertex-transitive, it must be $3$-cycle-regular. In particular, every dart in the fibre of $(bb)_r$ lies on a $3$-cycle and so, in $\Gamma$, there is a $3$-cycle $C$ through the edge $b_0b_r$. Then $\pi(C)$ is a $\lambda$-reduced closed walk of length $3$ that traces $(bb)_r$. It is straightforward to see that necessarily $\pi(C)=((bb)_r,(bb)_r,(bb)_r)$. Moreover, by Lemma \ref{cor:walks}, $0 \in \term(\pi(C))=\{3r\}$ and so 
\begin{align}
\label{eq:cong1} 3r \equiv 0 \quad (\mod 2m).
\end{align}
Now, every dart in the fibre of $(cc)_s$ must also lie on a $3$-cycle, and by an analogous argument,
\begin{align}
\label{eq:cong2} 3s \equiv 0 \quad (\mod m).
\end{align} 
Since $\Gamma$ is connected, by Lemma \ref{lem:ccv} we see that $\gcd(m,r,s) = 1$ and by (\ref{eq:cong1}) and (\ref{eq:cong2}), we see that the only possibility is that $m=3$, $r=2$ and $s=1$. One can readily verify that $\Gamma$ is isomorphic to the truncation of $K_{3,3}$, and thus is a vertex-transitive tricirculant of order $18$.     
\end{proof}

\begin{lemma}
If $\Gamma$ is a vertex-transitive $ccv$-cover of $\Delta_{7}$, then $\Gamma$ is a bicirculant of order $12$.
\end{lemma}

\begin{proof}
Suppose $\Gamma=\Cov(\Delta,\lambda,\iota,\zeta)$ is vertex-transitive where $(\Delta,\lambda,\iota,\zeta)$ is a $ccv$-extension of $\Delta_{7}$. Then $\iota(c)=\iota(a)=2m$ and $\iota(b) = \iota(d) = m$ for some $m \in \ZZ$. Let $W=((ad)_0,(da)_0,(ab)_0,(ba)_0)$ and observe that every reduced walk in $\mathcal{L}(W)$ is a $4$-cycle. In particular, $a_0b_0$ and $a_0d_0$ lie on a $4$-cycle, and thus, the vertex-transitivity of $\Gamma$ implies that there is a $4$-cycle $C$ through $c_0c_r$. It follows that $\pi(C)$ is a $\lambda$-reduced closed walk of length $4$ through the dart $(cc)_r$. It is plain to see that $\pi(C)=((cc)_r,(cc)_r,(cc)_r,(cc)_r)$ and thus $\term(\pi(C))=\{4r\}$. Then $4r \equiv 0$ ($\mod 2m)$ and $\gcd(m,r)=1$. Since $0 < r < m$, we see that $r = 1$ and $m=2$. Then $\Gamma$ is a cubic bicirculant of order $12$ and is in fact isomorphic to the Franklin graph. 
\end{proof}

\begin{lemma}
If $\Gamma$ is a vertex-transitive $ccv$-cover of $\Delta_{8}$, then $\Gamma$ is the triangular prism ${\rm GP}(3,1)$.
\end{lemma}

\begin{proof}
Let $\Gamma=\Cov(\Delta,\lambda,\iota,\zeta)$ be vertex-transitive where $(\Delta,\lambda,\iota,\zeta)$ is a $ccv$-extension of $\Delta_{8}$. Let $\iota(c)=m$ so that $\iota(a)=\iota(b)=2m$. Consider the walk $W=( (aa)_m,(ab)_0,(bb)_m,(ba)_0)$ and see that both $a_0a_m$ and $a_0b_0$ lie on a $4$-cycle in $\mathcal{L}(W)$. Since $\Gamma$ is vertex-transitive, then one of $d_0a_0$ or $d_0a_m$ must lie on a $4$-cycle $C$. Then $\pi(C)$ is a $\lambda$-reduced closed walk of length $4$ through the dart $(da)_{0}$. Clearly $\pi(C) = ((da)_0,(ab)_0,(bc)_0,(cd)_{-r})$ and $\term(\pi(C))=\{-r\}$. Then, by Lemma \ref{cor:walks} we have $r \equiv 0$ ($\mod m$). Since $\gcd(m,r)=1$ and $r < m$, we see that $m = 1$ and $r=0$. Then $\Gamma$ can be seen to be isomorphic to the triangular prism.
\end{proof}

\begin{lemma}
If $\Gamma$ is a $ccv$-cover of $\Delta_{9}$, then $\Gamma$ is not vertex-transitive.
\end{lemma}

\begin{proof}
Let $\Gamma=\Cov(\Delta,\lambda,\iota,\zeta)$ where $(\Delta,\lambda,\iota,\zeta)$ is a $ccv$-extension of $\Delta_{9}$. Observe that $\iota(a) = 2\cdot\iota(b) = 3\cdot\iota(d)$. Moreover, $\iota(b)$ is even as $b$ is incident to a semi-edge. Then, $\iota(a)$ is divisible by $12$. That is, for some $m \in \ZZ$ we have $\iota(a)=\iota(c)=12m$, $\iota(b) = 6m$, $\iota(d)=4m$ and the order of $\Gamma$ is $34m$. Suppose $\Gamma$ is vertex-transitive and consider the walk $$W=( (da)_0,(ab)_0,(ba)_0,(ad)_0,(da)_0,(ab)_0,(ba)_0,(ad)_0).$$ Observe that every dart beginning at $d_0$ lies on an $8$-cycle belonging to $\mathcal{L}(W)$. This implies the existence of a $\lambda$-reduced walk $W'$ of length $8$ through the dart $(bb)_{3m}$. Observe that then $W'$ must be one (or the inverse) of the following $6$ walks,  where $i \in \{-1,1\}$:
$$ W_{1,i} = ((bb)_{3m}, (ba)_0, (ac)_0, (cc)_{ir}, (cc)_{ir}, (cc)_{ir}, (ca)_0, (ab)_0),$$ 
$$ W_{2,i} = ((bb)_{3m}, (ba)_0, (ac)_0, (cc)_{ir}, (ca)_{0}, (ad)_{0}, (da)_0, (ab)_0),$$ 
$$ W_{3,i} = ((bb)_{3m}, (ba)_0, (ad)_0, (da)_{0}, (ac)_{0}, (cc)_{ir}, (ca)_0, (ab)_0).$$
Let $\mathcal{W}$ be the set containing the six walks $W_{j,i}$, $j \in \{1,2,3\}$ and $i \in \{-1,1\}$, along with their inverses. Denote by $\term(\mathcal{W})$ the union of endsets over the elements of $\mathcal{W}$. A tedious but straightforward computation shows that 
$$\term(\mathcal{W})=\{\pm 2m, \pm(m+r),\pm(m-r),\pm(3m+r),\pm(3m+3r)\}.$$
Then $z \equiv 0$ ($\mod 6m$) for some $z \in \term(\mathcal{W})$. This implies that $m=1$ and $r = km$ for some $k \in \{1,3,5\}$.  Then $\Gamma$ is one of three possible graphs of order $34$ (observe that $\Gamma$ is completely defined by the values of $r$ and $m$). One can check that in neither one of the three possible cases is $\Gamma$ vertex-transitive.
\end{proof}

\begin{figure}[h!]
\centering
\includegraphics[width=0.8\textwidth]{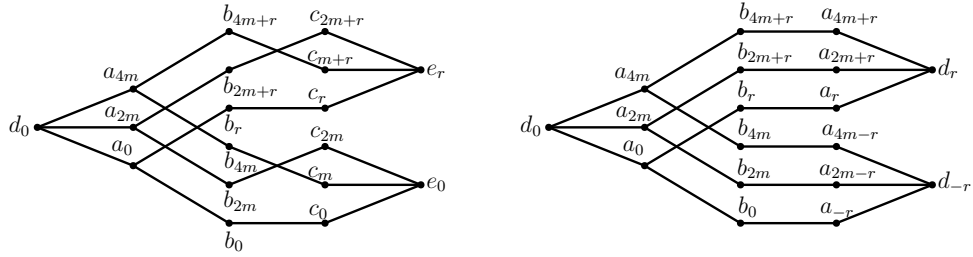}
\caption{The $8$ distinct $8$-cycles through the dart $d_0a_0$ in the proof of Lemma \ref{lem:delta19} }
\label{fig:delta10_1}
\end{figure}

\begin{lemma}
If $\Gamma$ is a $ccv$-cover of $\Delta_{10}$, then $\Gamma$ is not vertex-transitive.
\label{lem:delta19}
\end{lemma}

\begin{figure}[h!]
\centering
\includegraphics[width=0.8\textwidth]{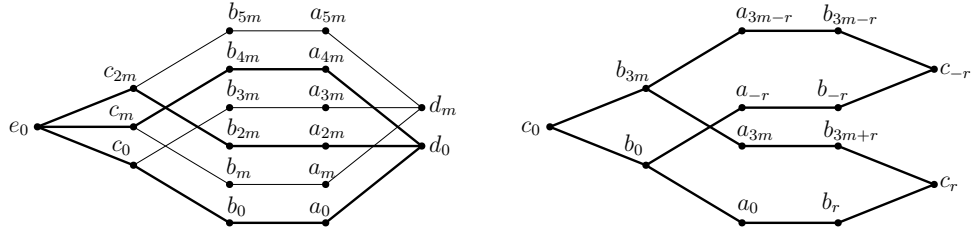}
\caption{Two subgraphs of $\Gamma$ in the proof of Lemma \ref{lem:delta19}. On the left, a subgraph containing all lifts of $W_1$ that visit $c_0$ ($8$-cycles through $c_0b_0$ shown in bold edges). On the right,  a subgraph containing all lifts of $W_3$ and $W_4$ that visit $c_0$.  }
\label{fig:delta10}
\end{figure}

\begin{proof}
Let $\Gamma=\Cov(\Delta,\lambda,\iota,\zeta)$ where $(\Delta,\lambda,\iota,\zeta)$ is a $ccv$-extension of $\Delta_{10}$. Suppose $\Gamma$ is vertex-transitive.
Since $\Delta_{10}$ has a $[1,3]$-edge, $\Gamma$ is arc-transitive.
 For some $m \in \ZZ$ we have $\iota(e)=m$, $\iota(c)=3m$, $\iota(a)=\iota(b)=6m$ and $\iota(d)=2m$. Observe that the order of $\Gamma$ is $18m$. As one can verify with the census of cubic vertex-transitive graphs \cite{census}, no $ccv$-cover of $\Delta_{10}$ is vertex-transitive if $m \in \{1,2\}$. Thus assume  that $m > 2$. Furthermore, $r \not \equiv 0$ ($\mod m$), for otherwise $m=1$ (since $\gcd(m,r)=1$).  Now, for $i,j \in \{0,1\}$ and $k \in \{2,4\}$ consider the walks in $\Gamma$: 
\begin{eqnarray*}
W_{i,j}&=&(d_0,a_{0},b_{ir},c_{ir},e_{ir},c_{(1+j)m + ir},b_{(4-2j)m + ir},a_{(4-2j)m},d_0),\\
W'_{i,k}&=&(d_0,a_0,b_{(1-i)r},a_{(1-2i)r},d_{(1-2i)r},a_{km+(1-2i)r},b_{km+(1-i)r},a_{km}). 
\end{eqnarray*}
%
Note that each of these $8$ walks is an $8$-cycle through $a_0d_0$ (see Figure \ref{fig:delta10_1}). Since $\Gamma$ is arc-transitive, then there must be $8$ distinct $8$-cycles through $c_0b_0$. Now, consider the walks
$$ W_1 = ((cb)_0, (ba)_{0}, (ad)_0, (da)_0, (ab)_0, (bc)_0, (cd)_0, (dc)_0),$$
$$ W_2 = ((cb)_0, (ba)_{-r}, (ad)_0, (da)_0, (ab)_{r}, (bc)_0, (cd)_0, (dc)_0),$$ 
$$ W_3 = ((cd)_0, (ba)_{0}, (ab)_{r}, (bc)_0, (cb)_0, (ba)_{-r}, (ab)_{0}, (bc)_0),$$
$$ W_4 = ((cd)_0, (ba)_{-r}, (ab)_{0}, (bc)_0, (cb)_0, (ba)_{0}, (ab)_{r}, (bc)_0)$$
(see Figure~\ref{fig:delta10}). Observe that for all $i \in \{1,2,3,4\}$, all lifts of $W_i$ that trace the dart $c_0b_0$ are $8$-cycles. There are exactly $6$ such cycles. It follows that there are an additional two $8$-cycles through $c_0b_0$ that do not project to any of the four walks $W_i$, $i \in \{1,2,3,4\}$. Each of these two cycles projects to a $\lambda$-reduced closed walk of length $8$ based at $c$ and visiting the dart $(cb)_0$. Moreover, such projections must be different than $W_i$ with $i \in \{1,2,3,4\}$. Let $\mathcal{W}$ be the set of all $\lambda$-reduced closed walks of length $8$ based at $c$ and visiting the dart $(cb)_0$, that are distinct from $W_i$ with $i \in \{1,2,3,4\}$. Let $\term (\mathcal{W})$ be the union of the endsets of all the elements of $\mathcal{W}$. A computer assisted calculation shows that $\term (\mathcal{W}) = \{\pm (2m + r), \pm(m+2r),\pm(m+r),\pm(r),\pm(2r),\pm(3r),\pm(2m+2r)\}$. Then $z \equiv 0$ ($\mod 3m$) for some $z \in \term \mathcal{W}$. This implies that $m=1$ and $r \in \{1,2,3,4,5\}$ or $m = 2$ and $r \in \{5,7\}$, a contradiction. We conclude that $\Gamma$ is not vertex-transitive.
\end{proof}

Up to this point we have shown that none of the labelled graphs $\Delta_i$ with $i \in \{1,\ldots,11\}$ admit a vertex-transitive $ccv$-covers with more than $20$ vertices. Since the graphs in Figure~\ref{fig:kcirc} all admit vertex-transitive $ccv$-covers with more than $20$ vertices, it follows that the set $\mQ$ consists of these eight graphs and possibly the graph $\Delta_{12}$ (it will be shown in Section \ref{sec:exception} that $\Delta_{12}$ does in fact admit infinitely many vertex-transitive $ccv$-covers and thus belongs to $\mQ$). Since the graphs in Figure~\ref{fig:kcirc} have at most $3$ vertices and all of their edges are of type $[1,1]$, we see that their generalised cyclic covers are $k$-multicirculants for some $k\le 3$.
The following proposition summarises the contents of this section.

\begin{proposition}
\label{prop:VT}
Let $\Gamma$ be a cubic vertex-transitive graph of order $n>20$ admitting an automorphism of order $n/3$ or greater. Then  either $\kappa(\Gamma) \in \{1,2,3\}$ or $\Gamma$ is a $ccv$-cover of $\Delta_{12}$.
\end{proposition} 

\subsection{The graph $\Delta_{12}$}
\label{sec:exception}

\begin{figure}[h!]
\centering
\includegraphics[width=0.2\textwidth]{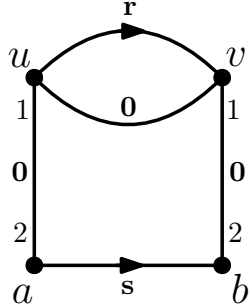}
\caption{The voltage assignment giving rise to the graph $\Gamma_{12}(m,r,s)$.}
\label{fig:type5}
\end{figure}

Let $m$ be a positive integer, and let $r$ and $s$ be two distinct elements of $\ZZ$. Let $\Delta_{12}(m,r,s)$ be the $ccv$-extension of $\Delta_{12}$ shown in Figure \ref{fig:type5}, where voltages are shown in bold characters next to each edge and $\iota(a) =m$. It follows from equality \ref{eq:ratio} that $\iota(u) = \iota(v) = 2m$ and $\iota(a) = \iota(b) = m$.
As the final step in the proof of Theorem~\ref{the:main}, we need to prove the following:

\begin{theorem}
\label{theo:delta11}
Let $m,r,s$, $m> 3$, be positive integers. If
 the cyclic generalised cover $\Gamma=\Gamma_{12}(m,r,s)$ arising from $\Delta_{12}(m,r,s)$ is connected and
 vertex-transitive, then $m$ is odd and $\Gamma \cong \SDW(m,3)$. Conversely, if
 $m$ is odd, then $\SDW(m,3)$ is isomorphic to $\Gamma_{12}(m,1,2)$, is connected and vertex-transitive and admits an automorphism of order $2m$.
\end{theorem}

\begin{proof}

 Recall that $\Gamma_{12}(m,r,s)$ admits an automorphism $\rho$ of order $2m$ whose orbits on vertices and darts are precisely the fibres of vertices and darts. 
Moreover, the automorphism $\rho^m$ fixes the $2m$ vertices in the fibres $\fib(a)$ and $\fib(b)$; in particular, $\Gamma_{12}(m,r,s)$ admits a non-trivial automorphism that
fixes one third of the vertices of the graph.
  
Suppose first that the girth of $\Gamma$ is less then $6$. Then by
Lemma~\ref{lem:girthregular} (see also \cite[Theorem 1.5]{signature})
 and the fact that $\Gamma$ has at least $24$ vertices, it follows that $\Gamma$ is isomorphic to the M\"obius ladder $\Lad(6m)$ or to the prism $\Prism(3m)$. However, every non-trivial automorphism of these graphs fixes at most $4$ vertices, yielding a contradiction.
  Hence, the girth of $\Gamma$ is at least $6$. Now observe that $\Gamma$ contains two $6$-cycles $C_1:=(u_0,a_0,u_m,v_m,b_m,v_0)$
and $C_2:=(u_0,a_0,u_m,v_{m+r},b_{m+r},v_r)$, implying that its girth is $6$. 

Let $(\alpha, \beta, \gamma)$, $\alpha\le \beta\le \gamma$, be the $6$-signature of $\Gamma$. Note that the edge $\{u_0,a_0\}$ lies on both of the cycles $C_1$ and $C_2$, while each of the remaining two edges incident with $u_0$ lies on precisely one of them. This shows that $\alpha, \beta \ge 1$ and $\gamma \ge 2$. Similarly, since both $C_1$ and $C_2$
pass through the edges $\{a_0,u_0\}$ and $a_0,u_m\}$, incident with $a_0$, we see that
$\beta, \gamma \ge 2$. Moreover, since $\alpha\ge 1$, there must a third $6$-cycle $C_3$ passing through the edge $\{a_0, b_s\}$. Since $C_3$ also passes through one of the edges
$\{a_0,u_0\}$, $\{a_0,u_m\}$, it follows that $\gamma \ge 3$. Finally, since $\Gamma$ admits
the automorphism $\rho^m$ (where $\rho$ is the canonical covering transformation) of order $2$ fixing a vertex and swapping two of its neighbours,
we see that two of the parameters $\alpha, \beta, \gamma$ must be equal.
By Lemma~\ref{lem:girthregular} (see also \cite[Theorem 1]{girth6VT}), it follows that
$\Gamma\cong \SDW(m,3)$. Finally, as mentioned in Remark~\ref{rem:g6} (and proved in 
\cite[Proposition 5]{girth6VT}), the automorphism group of $\SDW(m,3)$, $m\not = 3$, equals
$S_3 \times D_m$ and thus contains an element of order $2m$ if and only if $m$ is odd.

To conclude the proof of the theorem, assume that $m$ is odd, $m>3$, and consider
the mapping $\varphi \colon \V(\Gamma_{12}(m,1,2)) \to \V(\SDW(m,3))$ given by
\begin{eqnarray*}
\varphi(i,0,0)  & = & b_{i+1};    \\
\varphi(i,0,1) & =  & a_i;          \\
\varphi(i,1,0)  & = & 
     \left\{ \begin{array}{ll} 
        u_{i+m} & \hbox{ if } i \hbox{ is even}; \\   u_{i} & \hbox{ if } i \hbox{ is odd}; 
      \end{array}\right.        \\
\varphi(i,1,1) & =  & 
     \left\{ \begin{array}{ll} 
        v_{i+1} & \hbox{ if } i \hbox{ is even}; \\   v_{i+m+1} & \hbox{ if } i \hbox{ is odd}; 
      \end{array}\right.       \\
\varphi(i,2,0)  & = & 
   \left\{ \begin{array}{ll} 
        u_{i} & \hbox{ if } i \hbox{ is even;} \\   u_{i+m} & \hbox{ if } i \hbox{ is odd;} 
      \end{array}\right.        \\
\varphi(i,2,1) & =  & 
     \left\{ \begin{array}{ll} 
        v_{i+m+1} & \hbox{ if } i \hbox{ is even;} \\   v_{i+1} & \hbox{ if } i \hbox{ is odd;} 
      \end{array}\right.        
\end{eqnarray*}
for all $i \in \{0,1, \ldots, m-1\}$, where the indices at $b_j$'s and $a_j$'s are computed modulo $m$, while those at $u_j$'s and $v_j$'s are computed modulo $2m$. 
The isomorphism $\varphi$ is depicted in Figure~\ref{fig:wreath} where each vertex of $\SDW(m,3)$
is labeled with its $\varphi$-image. 
\begin{figure}[h!]
\centering
\includegraphics[width=0.7\textwidth]{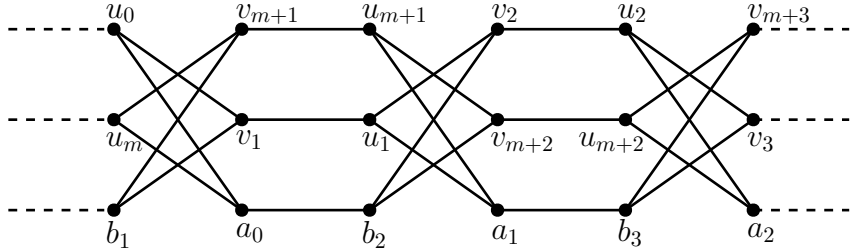}
\caption{A section of $\Gamma_{12}(m,1,2)$ with $m \geq 3$}
\label{fig:wreath}
\end{figure}
It is obvious that $\varphi$ is a graph isomorphism.
\end{proof}

By \cite[Proposition 5]{girth6VT} (see also Remark~\ref{rem:g6}), the automorphism group of 
$\SDW(m,3)$ is isomorphic to $S_3\times D_m$, unless $m=3$, in which case $\SDW(m,3)$ is the unique cubic arc-transitive graph on
$18$ vertices, also called the Pappus graph. From this, it is easy to deduce the parameter
$\kappa$ and $\relmeo$ for the split wreath graph $\SDW(m,3)$:

\begin{lemma}
\label{lem:deltanomulti}
Let $\Gamma = \SDW(m,3)$ for some integer $m\ge 3$. Then one of the following holds:
\begin{itemize}
\item $m$ is not divisible by $3$ and 
$\relmeo(\Gamma) = \kappa(\Gamma) = 2$;
\item
 $m \equiv 0$  $(\mod 6)$, and $\relmeo(\Gamma) = \kappa(\Gamma) = 6$;
\item
$m \equiv 3$  $(\mod 6)$, $m\not = 3$, and
 $\relmeo(\SDW(m,3)) = 3$ while $\kappa(\SDW(m,3)) = 6$;
\item 
$m=3$, $\Gamma$ is isomorphic to the Pappus graph 
and $\kappa(\Gamma) = 3$ while $\relmeo(\Gamma) = 3/2$.
\end{itemize}
\end{lemma}

\color{black}

\section{Proof of Theorem~\ref{the:main}}
\label{sec:last}
Using the results proved in the previous sections, it is now not difficult to prove
Theorem~\ref{the:main}. Before we proceed to the proof, we will need the lemma below, characterizing those cyclic Haar graphs that are circulant graphs. 

\begin{lemma}
\label{lem:haar}
For an integer $m \geq 5$, a connected cyclic Haar graph $\rH(m,x,y)$ is a circulant if and only if $m$ is odd and $\{x,y\} = \{a,2a\}$ or $\{x,y\} = \{a,-a\}$ for some $a \in \ZZ$ such that $\gcd(m,a)=1$.
\end{lemma}

\begin{proof}
Suppose that $\Lambda=\rH(m,x,y)$, $\gcd(m,x,y) = 1$, is a circulant. For $i\in \ZZ_m$, we let $u_i$ and $v_i$ denote the vertices $(i,0)$ and $(i,1)$, respectively, and recall that the neighbours of $u_i$ in $\Lambda$ are $v_i, v_{i+x}$ and $v_{i+y}$.

Since $\Lambda$ is a circulant, it is isomorphic to a prism or to a M\"{o}bius ladder. In both cases, the girth of $\Gamma$ is $4$
and at every vertex there is an edge belonging to two $4$-cycles, none of these 
$4$-cycles sharing a $2$-path (recall that we are assuming that $m \geq 5$ and thus the order of $\Lambda$
is at least $10$). Observe also that exactly one of the graphs
$\Lambda$, $\rH(m,-x, y-x)$ and $\rH(m, x-y, -y)$, isomorphic to $\Lambda$, is such that
the edge incident with the vertex $(0,0)$ belonging to two $4$-cycles is $e_0=\{u_0, v_0\}$.
Denote this graph $\rH(m,a,b)$ (that is, $(a,b)$ is one of $(x,y), (-x,y-x)$ or $(x-y,-y)$,
or equivalently, $(x,y)$ is one of $(a,b)$, $(-a,b-a)$ or $(a-b,-b)$).

One of the two $4$-cycles passing through $e_0$ passes also through the edge $e_1=\{u_0,v_a\}$, while the other passes through $e_2=\{u_0,v_b\}$. 
Since the neighbours of $v_0$ are $u_0, u_{-a}$ and $u_{-b}$,
for these two edges to form
a $4$-cycle together with the edge $e_0$, 
either $u_{-a}$ is adjacent to $v_a$ and $u_{-b}$ to $v_{b}$, or
$u_{-b}$ is adjacent to $v_a$ and $u_{-a}$ to $v_{b}$.
In the first case, both $2a$ and $2b$ belong to the set $\{0,a,b\}$,
while in the second case, $a+b \in \{0,a,b\}$, which forces $a+b = 0$
(since $a\not = 0 \not = b$).

Let us consider the first case.
Observe that $2a \not = a$ and $2b \not = b$.
If  $2a=0$, then $m$ is even, $a=m/2$ and $2b=a$. Since $\gcd(m,a,b) = 1$, this forces
$m=4$, contradicting our assumption that $\Lambda$ has at least $10$ vertices.
A similar contradiction is obtained if $2b=0$. This leaves us with the possibility that
$2a=b$ and $2b =a$. But then $a + b =0$, and we find ourselves in the second case.

Therefore, the second case occurs, that is, $a+b = 0$. Since $\gcd(m,a,b) = 1$,
this implies that $\gcd(m,a) = 1$ and $b=-a$. But then $(x,y)$ is one of the pairs
$(a,-a)$, $(-a,-2a)$ or $(2a,a)$.

To summarise, if a connected cubic cyclic Haar graph $\rH(m,x,y)$ with $m\ge 5$ is
a circulant then $\{x,y\} = \{a, 2a\}$ or $\{x,y\} = \{a, - a\}$ 
for some $a \in \ZZ_m$ with $\gcd(m,a)=1$. Furthermore, the graphs $\rH(m,a,2a)$ and $\rH(m,a,-a)$ are isomorphic to the prism $\Prism(m)$ whenever $m$ is even, and thus cannot be circulants. It follows that $m$ must be odd, as required.

For the converse, let $m,a \in \ZZ$ be such that $m \geq 5$, $m$ is odd and $\gcd(m,a)=1$. If $\Lambda = \rH(m,a,-a)$ then the mapping given by $u_i \mapsto v_{i+r}$ and $v_i \mapsto u_{i+r}$ is a circulant automorphism of $\Lambda$. On the other hand if $\Lambda = \rH(m,a,2a)$ then the mapping given by $u_i \mapsto v_i$ and $v_i \mapsto u_{i-2r}$ is a circulant automorphism of $\Lambda$.
\end{proof}

We are now ready to prove Theorem \ref{the:main}. As mentioned in Section~\ref{sec:strategy}, the validity of the
theorem for graphs on at most $20$ vertices can easily be checked by consulting the census of cubic vertex-transitive graphs \cite{census}.

Now, as in the statement of Theorem~\ref{the:main}, let $\Gamma$ be a cubic vertex-transitive graph of order $n$, $n>20$, admitting an automorphism of order at least $\frac{n}{3}$. We need to show that then one of the claims (1) -- (4) of Theorem~\ref{the:main} holds.

Observe first that by combining Proposition~\ref{prop:VT} with Theorem~\ref{theo:delta11}, 
$\kappa(\Gamma) \in \{1,2,3\}$ or $\Gamma$ is isomorphic to
$\SDW(m,3)$ with $m$ odd (and $m\not = 3$ since we are assuming that $\Gamma$ has
more than $20$ vertices).
 In the latter case, by Lemma~\ref{lem:deltanomulti}, we see that
$\kappa(\Gamma) \le 3$ or $m \equiv 3\> (\mod 6)$ (in which case $\kappa(\Gamma) = 6$).
 %
To summarise, either $\kappa(\Gamma) \in \{1,2,3\}$ or claim (4) of Theorem~\ref{the:main} holds.

Suppose now that $\kappa(\Gamma) \in \{1,2,3\}$.
If $\Gamma$ is a circulant, then, by definition, $\Gamma \equiv \Cay(\ZZ_{2m};S)$, where
$S=\{r,-r,m\}$ for some $r\in \ZZ_{2m}$, $r\not = 0$. By connectivity of $\Gamma$, we may assume that $\gcd(r,m) = 1$, implying that there exists $s\in \ZZ_{2m}$, $\gcd(s,2m) =1$, such that
$rs = 1$, or $r$ is even, $m$ is odd and $rs=2$ in $\ZZ_{2m}$. Note that $\Gamma$ is then isomorphic
to $\Cay(\ZZ_{2m};\{1,-1,m\}) \cong \Lad(2m)$ (in the first case) or to $\Cay(\ZZ_{2m};\{2,-2,m\}) \cong \Prism(m)$ (in the second case). In particular, claim (1) of
Theorem~\ref{the:main} holds in this case.

If $\Gamma$ is a bicirculant, then it can be deduced from
\cite[Propositions 3 and 4]{bic} and \cite[Theorem 7 and Corollay 8]{Igraphs}
(see also \cite[Theorem 1.1 and Remark 1.2]{cycliccovers}) that
$\Gamma$ is  isomorphic either to a
prism $\Prism(m)$, 
or to a M\"obius ladder $\Lad(n)$, 
or to one of the vertex-transitive
generalised Petersen graphs $\GP(m,r)$ with $2\le r < m/2$, $r^2\equiv \pm 1\> (\mod m)$,
or to a cyclic Haar graph $\rH(m,r,s)$ with $\gcd(m,r,s) = 1$.

If $\Gamma$ is a prism or a M\"obius ladder, then it is  not 
a circulant only if it is isomorphic to a prism $\Prism(m)$ with $m$ even (and then
claim (2a) of Theorem~\ref{the:main} holds).
Further, since the girth of a generalised Petersen graph $\GP(m,r)$ is at least $5$ unless $m\le 4$ or $r=1$, we see that a generalised Petersen graph $\GP(m,r)$
is a circulant if and only if $r=1$ and $m$ is odd (in which case it is a prism).
In particular, if $\Gamma$ is a generalised Petersen graph which is not a circulant,
then claim (2b)  of Theorem~\ref{the:main} holds.

If $\Gamma$ is 
a connected cyclic Haar graph, that is $\Gamma=\rH(m,r,s)$ with $m> 10$ and
$\gcd(m,r,s) = 1$, then by applying the permutation $\varphi_{\alpha,0}$ for an appropriate $\alpha$, we may assume that $r$ divides $m$ (where $r$ is represented as a positive integer smaller than $m$),
and thus that $\gcd(r,s) = 1$. Under this assumption,
by Lemma \ref{lem:haar}
 $\Gamma$ is then a circulant if and only if $m$ is odd and
$\{r,s\} = \{1,  2\}$ or
$\{r,s\}=(1,m-1)$.
We have thus shown that the claim (2c) of Theorem~\ref{the:main} holds in this case.

If $\Gamma$ is a tricirculant but not a bicirculant, then by 
\cite[Theorems 1.1, 4.3 and 5.3]{tricirc}, 
$\Gamma$ is either
the Tutte's $8$-cage (on $30$ vertices),
 the truncated tetrahedron (on $12$ vertices) or
isomorphic to one of the graphs
$\X(m)$ or $\Y(m)$ with $m \equiv 3\> (\mod 6)$. Note that claims (3) of
Theorem~\ref{the:main} holds in this case.

For the converse, it is clear that the circulants, bicirculants and tricirculant appearing in parts (1), (2) and (3) of Theorem~\ref{the:main} all admit an automorphism of order at least one third of the order of the graph. The graph $\SDW(m,3)$ with $m \equiv 3\> (\mod 6)$ has order $6m$ and admits, by construction, an automorphism with two orbits of size $2m$, namely, the canonical covering transformation $\rho$. This completes the proof of Theorem~\ref{the:main}.

\section{Relationship between $\kappa(\Gamma)$ and $\eta(\Gamma)$ and open problems}
\label{sec:problem}

Let us conclude the paper with a discussion on the interplay
between  the parameters $\kappa(\Gamma)$ and $\relmeo(\Gamma)$.
First, note that none of the functions $\relmeo(\Gamma)$ and
$\kappa(\Gamma)$ can be bounded above by a constant, even when restricted to
the class of cubic vertex-transitive graphs. Namely, if there were a constant $C$ such
that $\relmeo(\Gamma) \le C$ holds for all graphs $\Gamma$,
then $\max \{o(g) : g \in \Aut(\Gamma)\} \ge |\V(\Gamma)|/C$, or in other words,
the parameter $\meo(\Gamma):=\max \{o(g) : g \in \Aut(\Gamma)\}$ is bounded 
below by a linear function of $|\V(\Gamma)|$. However, it was shown in \cite{regorbs} that there exists an infinite family $\{\Gamma_i\}_{i\in\NN}$ of cubic vertex-transitive graphs with $\meo(\Gamma_i)$ growing slower than any logarithmic function of $|\V(\Gamma_i)|$. This implies that
there is no constant $C$ such that $\relmeo(\Gamma) \le C$ for all cubic vertex-transitive
graph $\Gamma$. Since $\relmeo(\Gamma)  \le \kappa(\Gamma)$ holds for all $\Gamma$,
this implies that there is no upper bound on $\kappa(\Gamma)$.

Let us mention at this point
that a somewhat similar problem was raised in \cite{CamSSpi},
where it was conjectured that the value
$$ m(n):=\min \{ \frac{|\V(\Gamma)|}{\kappa(\Gamma)} : \Gamma \hbox{ a cubic vertex-transitive graph on at most } n \hbox{ vertices} \}$$
tends to $\infty$ as $n$ grows to $\infty$, or in other words, 
for every integer $n$, there is a constant $c_n$ such that
every cubic vertex-transitive graph on more than $n$ vertices admits
a semiregular element of order at least $c_n$.
However, in \cite{SpiCubicPoly}, an infinite family of cubic vertex-transitive graphs was
constructed in which every semiregular automorphism has order at most $6$, proving this
conjecture wrong in general (but proved to be correct when restricted to arc-transitive or Cayley graphs; see also \cite{SPX}).

While $\kappa(\Gamma)$ can take an arbitrary large value,
a question arises whether the parameter $\kappa(\Gamma)$ can be bounded
above in terms of $\relmeo(\Gamma)$ where $\Gamma$ ranges through the class of cubic vertex-transitive graphs. 
More precisely, we are interested in the following:

\begin{question}
\label{q:kr}
For which positive integers $r$ does there exist an integer $k_r$ such that
$\kappa(\Gamma) \leq k_r$ for all but finitely many cubic vertex-transitive graphs satisfying $\relmeo(\Gamma) \leq r$.
\end{question}

If for some $r$ an integer $k_r$ as above exists, then let $f(r)$ be the smallest such integer $k_r$.
A direct consequence of Theorem~\ref{the:main} is that $f(r)$ can be defined at least for $r\in \{1,2,3\}$ and that $f(1) = 1$, $f(2) =2$ and $f(3)=6$. 
Let us conclude this paper by posing the following:

\begin{question}
If the $f(r)$ is defined for all $r$, what is its asymptotic behaviour as $r\to \infty$? Is $f(r)$ unbounded?
If so, can it be bounded by a linear and/or polynomial function of $r$?
\end{question}

\end{document}